\documentclass{amsart}

\usepackage{graphicx}
 
\def\R{{\mathbb R}}
\def\N{{\mathbb N}}

\def\Z{{\mathbb Z}}
\def\1{{1\!\!\!1}}
\def\E{{\mathbb E}}
\def\P{{\mathbb P}}

\def\cal{\mathcal}

\def\ol{\overline}

\def\mes{{\rm{mes}}}

\def\eps{\varepsilon}
\newtheorem{theorem}{Theorem}
\newtheorem{prop}{Proposition}[section]
\newtheorem{cor}{Corollary}[section]
\newtheorem{defi}{Definition}[section]
\newtheorem{lemma}{Lemma}[section]

\title{On the spectrum of Markov semigroups \\ via sample path 
large deviations.}
\author{Irina Ignatiouk-Robert}
\address{
{Universit\'e de Cergy-Pontoise,}
{D\'epartement de math\'ematiques,}
{2, Avenue Adolphe Chauvin,}
{95302 Cergy-Pontoise Cedex,}
{France}}
\date{May, 20 2002}
\email{Irina.Ignatiouk@math.u-cergy.fr}
 
\begin{document}
\begin{abstract}
The essential spectral radius of a sub-Markovian process is defined
as the infimum of the spectral radiuses of all local perturbations of
the process. When
the family of rescaled processes satisfies sample path large deviation
principle, the spectral radius and the essential spectral radius are expressed in terms of
the rate function.  The paper is motivated by applications to reflected diffusions and
jump Markov processes describing stochastic networks for which the 
sample path large deviation principle has been established and the rate function has been 
identified while essential spectral radius has not been calculated. 
\end{abstract}
\maketitle
 
 \section{Introduction}\label{sec-Int}
 
For a sub-Markovian process $(X(t))$ on a locally compact set $E$ endowed with a
non-negative Radon measure $m$, {\it spectral radius} $r^*$
is defined as the infimum  over all those $r>0$ for which the resolvent function  
\[
R_r\1_W(x) = \int_0^\infty r^{-t} \P_x(X(t)\in W) \,dt 
\]
is  $m$-integrable on compact subsets of $E$ for every compact
set $W\subset E$. Under some general assumptions, the quantity $r^*$ can be described
in several ways: 
\begin{itemize}
\item[(i)] 
For an irreducible discrete time Markov chain $(X(n))$ on a countable set $E$,
$1/r*$ is a common radius of convergence of the  series
\[
\sum_{n=1}^\infty z^n p(n,x,y), \quad x,y\in E,
\]
where $\{p(n,x,y), x,y\in E\}$ are the transition probabilities of
$(X(n))$ (see Seneta~\cite{Seneta} and Vere-Jones~\cite{Vere-Jones:1}).
\item[(ii)]  
\[
\log r^* ~=~ \sup_{W,V} \limsup_{t\to +\infty} \frac{1}{t} \log m(\1_W P^t\1_V) 
\] 
where the supremum is taken over all compact subsets $V,W\subset E$ and $\{P^t,t>0\}$ is
the sub-Markovian  semi-group associated to $(X(t))$. 
\item[(iii)] $r^*$ is  the infimum over all  $r>0$ for
which there is a positive measurable function $f$  which is $m$-integrable on compact
subsets of $E$ and such that  $P^t f \leq r^t f$ for all $t>0$.  A
function $f$ satisfying  the inequality $P^t f \leq  r^t f$ is usually
called $r$-superharmonic. A dual description of $r^*$  can be given
by using $r$-superharmonic Radon measures on $E$ (see
Seneta~\cite{Seneta},            Stroock~\cite{Stroock}            and
Vere-Jones~\cite{Vere-Jones:1}  for example). 
\item[(iv)] When  the sub-Markovian
semi-group  $\{P^t,  t>0\}$ is generated by a symmetric linear operator
$A$  in  $L^2(E, m)$, the value $-\log   r^*$   is  the   bottom of   the
$L^2$-spectrum  of $A$ (see Stroock~\cite{Stroock} and also LimingWu~\cite{LimingWu} for a 
similar result for discrete time Markov chains).
\item[(v)] 
\[
r^* = \sup_K \left\{ r > 0 : \E_{\bf \, \cdot}( r^{-\tau_K} ) \in L^1(K,m) \right\} 
\]
where the supremum is taken over all compact subsets $K\subset E$ and for every compact
set $K$, $\tau_K$ denotes the first exit time from the set $K$ (see Stroock~\cite{Stroock}). 
\end{itemize}
The last description of the quantity $r^*$ shows that $r^*$ provides the rate at which the
process $(X(t))$ leaves compact sets. This quantity is of interest for 
transient Markov processes, because it  shows  how fast the process goes to infinity. 
 
Because  of the  first property,  the value  $1/r^*$ is  usually called
{\em convergence  parameter}   of  $(X(t))$\footnote{For  the
  definition  of convergence parameter  for general state space Markov  chains,
  see  Nummelin~\cite{Nummelin1}}.  In the present paper, Woess's terminology is 
used :  we call the quantity $r^*$ spectral radius of the process $(X(t))$, see  Woess's
book~\cite{Woess}.  While
this terminology  may be  misleading (in  a non-symmetrical case,  the correspondence  between
the  value  $r^*$ and  operator properties  of $\{P^t,t>0\}$ is more
intricate, see Vere-Jones~\cite{Vere-Jones:2}),  it stresses  the  importance of  the
correspondence  between the  value $r^*$  and the spectral  radius in
$L^2(E, m)$  of  a   symmetric  case.

The value $r^*$ is clearly sensitive to changing the transition probabilities on 
compacts.  The quantity 
$r^*_e$  is defined as the infimum of  $r^*$ over all such changes :   
\[
r^*_e = \inf_K r^*_K 
\]
 $r^*_K$  denotes here spectral radius of  the sub-Markovian process $(X(t))$
 killed on the set $K$ and the infimum is taken over all compact subsets $K\subset E$.

For symmetric Markov processes, by Perssons principle (see Grillo~\cite{Grillo:01} for
symmetric diffusions and Liming
Wu~\cite{LimingWu} for symmetric Markov chains),  the quantity $r^*_e$ is related to
the $L^2$-essential spectral radius of 
the corresponding Markov semi-group.  It is of interest for recurrent Markov processes :
given a compact set $K\subset E$, let $\tau(K)$ denote the first hitting time of $K$, then 
under some general assumptions,  the number $r^*_e$  equals  the infimum
over  all those $r>0$ for which the function  
\[
R_{K,r}\1(x) = \int_0^\infty r^{-t} \P_x(\tau(K) > t)  \,dt 
\]
is  $m$-integrable on compact subsets of $E\setminus K$ (see Proposition~\ref{pr3-3}
below). The quantity
$r^*_e$ provides therefore the rate at which the process returns to compacts. 
 
For some positive recurrent countable  Markov chains, the quantity $r^*_e$ gives an
accurate bound to the rate of convergence to equilibrium : Malyshev  and
Spieksma~\cite{Malyshev-Spieksma} have shown that this is the  best geometric
convergence rate  when the transitions  of  the Markov chain are changed on
a  finite  subset of $E$.   For more details concerning a relationship between the quantity
$r^*_e$ and the rate of convergence to equilibrium see 
Liming Wu~\cite{LimingWu}).

Unfortunately, in practice, outside of some particular examples, an explicit
representation of $r^*_e$ is very difficult to obtain. 
 
In the present paper, we consider a 
Markov process $(X(t))$ on $\R^d$ for  
which the family of rescaled processes $(Z_a(t)=X(at)/a, \; t\in[0,1])$ satisfies sample path
large deviation principle with a good rate function $I_{[0,1]}$. The quantities
$r^*$ and $r^*_e$ are represented  in terms of the rate function : we show that 
\begin{equation}\label{e1-I}
\log r^* = - \inf_{\phi : \phi(0) = \phi(1)} I_{[0,1]}(\phi) \quad \text{ and }  \quad \log
r^*_e = - \inf_{\substack{\phi(0)=\phi(1),  \\
\phi(t)\not= 0, \; \forall \; 0<t<1} }  I_{[0,1]}(\phi)
\end{equation}
where the first infimum is taken over all continuous functions $\phi$ with
$\phi(0)=\phi(1)$ and the second infimum is taken over all continuous functions $\phi$ with
$\phi(0)=\phi(1)$ such that $\phi(t)\not=0$ for all $0<t<1$. 
 
The first result in this domain was obtained  by  Malyshev and  
Spieksma~\cite{Malyshev-Spieksma} for  discrete time
partially homogeneous random walks  in $\N$ and in $\Z^2_+$.  
Unfortunately, their proofs use particular properties of the processes and it is
not usually possible to extend them to a more general situation (see section~\ref{sec1} for more
details). 
 
Our results are motivated by applications to reflected diffusions considered by Varadhan
and Williams~\cite{Varadhan-Williams} and  jump Markov processes
describing stochastic networks. For these processes,  the 
sample paths large deviation principle has been established and an explicit representation of
the  corresponding rate function has been obtained
(see~\cite{A-D,Delcoigne:01,D-E,Ignatiouk:01,Ignatiouk:04, S-W} for 
example) while  the essential spectral radius $r^*_e$ has not been identified. 
 
An example of Jackson networks illustrates our results. Using Relation~(\ref{e1-I}) we
obtain an explicit representation of the quantities $r^*$ and $r^*_e$ for two-dimensional
Jackson networks. In the forthcoming paper, we apply our results for reflected diffusions
in $\R_+^2$. 
 
 \section{The main results}\label{sec1}
We consider  a strong  Markov process $(X(t))$ on $E\subset\R^d$ whose sample
paths are right continuous with left limits. The set $E$ is endowed by a non-negative
Radon measure $m$. We  assume that the set $E$ is closed and unbounded and that
$m(O) \not= 0$ for every open non-empty subset $O\subset E$. Given a closed set $V\subset
E$, $\tau(V)$ denotes the hitting   
time of $V$: $$\tau(V)=\inf\{t > 0 : X(t)\in V\},$$ by convention $\tau(\emptyset) = +\infty$.
It is assumed that for  every real bounded measurable function $\varphi$ on $E$, the mapping 
$$(t,x)\to \E_x\bigl(\varphi(X(t)) ; \; \tau(V)>t \bigr)$$
is ${\cal B}(\R_+)\times{\cal B}(E)$-measurable from $\R_+\times E$ to
$\R$.

\begin{defi} {\bf Spectral radius} $r^*$ is the infimum of all those $r>0$ for which
  the  resolvent function 
\[
R_r\1_W(x) = \int_0^\infty r^{-t} \P_x(X(t)\in W )  \,dt 
\]
is $m$-integrable on compact subsets of $E$ for every compact subset $W\subset E$.   
\end{defi}
\begin{defi} {\bf Essential spectral radius} $r_e^*$ is the infimum of all those $r>0$ for which
there is a compact set $K\subset E$ such that the truncated resolvent function 
\[
R_{K,r}\1_W(x) = \int_0^\infty r^{-t} \P_x(X(t)\in W, \, \tau(K) > t)  \,dt 
\]
is $m$-integrable on compact subsets of $E\setminus K$ for every compact subset $W\subset
E\setminus K$. 
\end{defi}
 
The Markov process $(X(t))$ is assumed to satisfy the following large deviation conditions. 
 
\medskip 
\noindent
{\bf Assumption~(A) : Large deviations.} {\em Let ${\cal E}$ be the
set of all possible limits $\lim_{a\to\infty} x_a/a$ with $x_a\in E$, and let $D([0,T],{\cal E})$
denote the Skorohod space of all functions $\phi$ from $[0,T]$ to ${\cal E}$  which are right
continuous and have left limits. We assume that 
\begin{itemize}
\item[$(a_0)$] $E\subset{\cal E}\not=\{0\}$ and the set ${\cal E}\setminus\{0\}$ is convex;
\item[$(a_1)$] for every $T>0$, the family of
rescaled processes $$(Z_a(t), \; t\in[0,T])\stackrel{\text{def.}}{=} (X(at)/a, \; t\in[0,T])$$ satisfies 
sample
path large deviation principle in $D([0,T],{\cal E})$ with a good rate functions
$I_{[0,T]}$ (see section~\ref{sec4} for a precise definition);   
\item[$(a_2)$] the rate function $I_{[0,T]}$ has an integral form~: there is a local rate function
$L:{\cal E}\times\R^d\to\R_+$ such that 
\[
I_{[0,T]}(\phi) = \int_0^T L(\phi(t),\dot\phi(t))  \,dt 
\]
if the function $\phi : [0,1]\to {\cal E}$ is absolutely continuous, and $I_{[0,1]}(\phi)
= +\infty$ otherwise.
\item[$(a_3)$] there are two convex functions $l_1$ and $l_2$
  on $\R^d$   such that  
\begin{itemize}
\item[--]   $0\leq l_1(v) \leq L(x,v) \leq l_2(v)$ for all $x\in{\cal E}$ and for
  all $v\in\R^d$, 
\item[--] the function $l_2$ is finite in a neighborhood of zero 
\item[--] and $$
\lim_{n\to\infty} \inf_{|v|\geq n} l_1(v)/|v| > 0.$$
\end{itemize}
\end{itemize}}
For $x,y\in{\cal E}$ and $t>0$, we denote by $I(t,x,y)$ the infimum of the
rate function $I_{[0,t]}(\phi)$ over all continuous functions $\phi : [0,t] \to {\cal E}$
with $\phi(0)=x$ and $\phi(t)=y$,  $\hat{I}(t,x,y)$ denotes
the infimum of $I_{[0,t]}(\phi)$ over all continuous  functions $\phi : [0,t] \to {\cal E}$ for which
$\phi(0)=x$, $\phi(t)=y$ and the set $\{s\in[0,t]~: \phi(s)=0 \}$
has Lebesgue measure zero, $I_0(t,x,y)$ denotes
the infimum of $I_{[0,t]}(\phi)$ over all continuous  functions $\phi : [0,t] \to {\cal E}$ such that
$\phi(0)=x$, $\phi(t)=y$  and $\phi(s)\not=0$ for all $0<s<t$. 
$\ol{0} : \R_+\to \R^d$ denotes the constant function $\ol{0}(t)\equiv 0\in\R^d$. The
quantities $I^*$ and  $I^*_0$ are defined by   
\[
I^* = I(1,0,0) \quad \quad \text{and} \quad \quad I^*_0 = \hat{I}(1,0,0). 
\]
 
Using classical large deviation techniques we obtain the following result.
 
\begin{theorem}\label{th1} Under the hypotheses (A), for any $x,y\in{\cal E}$,
\begin{equation}\label{e1-1}
\log r^* ~=~ -I^* ~=~ - I_{[0,1]}(\ol{0})  ~=~ - \limsup_{t\to +\infty} I(t,x,y)/t  ~=~ -
\inf_{\phi(0)=\phi(1)} I_{[0,1]}(\phi)  
\end{equation}
where the infimum is taken over all continuous functions $\phi : [0,1]\to {\cal E}$ with
$\phi(0)=\phi(1)$.
\end{theorem}
 
Our main technical result is the following theorem. 
 
\begin{theorem}\label{th2} Under the hypotheses (A), for any $x,y\in{\cal E}$, 
\begin{align}
\log r^*_e ~=~ -I^*_0 &=~ - \limsup_{t\to +\infty} I_0(t,x,y)/t ~=~ - \limsup_{t\to +\infty}
\hat{I}(t,x,y)/t \label{e2-1} \\ &=~ - \inf_{\substack{\phi(0)=\phi(1),  \\
\phi(t)\not= 0, \; \forall \; 0<t<1} }  I_{[0,1]}(\phi) \nonumber
\end{align}
where the infimum is taken over all continuous functions $\phi : [0,1]\to {\cal E}$ with
$\phi(0)=\phi(1)$ such that $\phi(t)\not=0$ for all $0<t<1$.
\end{theorem} 
 
Theorem~\ref{th2} extends the result obtained earlier by Malyshev and
 Spieksma~\cite{Malyshev-Spieksma} for
 discrete time homogeneous random walks   in $\N$ and in
$\Z^2_+$. The main difficulty
consists here in the proof of the upper bound $\log r^*_e \leq - I^*_0$.
To get this inequality, one has to analyze the behavior of the rescaled processes $(Z_n(t))$
in a neighborhood of infinity where truncations on compact sets are not
sufficient. In a such situation, Freidlin-Wentzel method can not be applied.  
 
The proof of the upper bound $\log r^*_e \leq - I^*_0$ proposed by Malyshev
and Spieksma in \cite{Malyshev-Spieksma} uses particular properties of the
processes : they considered discrete time random walks with uniformly bounded jumps
for which the sets 
\[
K_x = \{ y: I(1,x,x+y)< +\infty\}
\] 
are compact and bounded uniformly in $x$. For continuous time Markov processes, the sets
$K_x$ are usually 
not bounded. Moreover, their method required that for any $\eps>0$ there exist
$\delta>0$ and $n_\eps>0$ such that for all $n > n_\eps$,
\[
 \sup_{|x'-x|< n\delta}
\frac{1}{n} \log \P_{x'}\left( |X(n) - n y| < \delta n \right) \; \leq
    \; - I(1,x,y) + \eps
\]
uniformly in $(x,y)$ on the set of all $(x,y)$ for which the quantity $I(1,x,y)$  is
finite. Such an uniform convergence is very difficult to 
check in practice and  is sometimes wrong~:  this implies the uniform continuity of the mapping
\[
(x,y) \to I(1,x,y) = \inf_{\phi(0)=x, \; \phi(1)=y} I_{[0,1]}(\phi).
\]
For  the standard  Brownian  motion  in $\R$,  this infimum  is
achieved  by  the function $\phi(t)=x  +   t(y-x)$  and  equals
$(y-x)^2/2$. The function $(y-x)^2/2$ is not uniformly continuous on $\R^2$.
 
In  our paper, we prove the inequality $\log r^*_e \leq - I^*_0$ by using a method of
statistical physics called cluster expansions,
see  for  example Malyshev  and  Minlos~\cite{Malyshev-Minlos}  or
Rivasseau~\cite{Rivasseau:01}. In the present setting, this method consists in
bounding the quantity $\log r^*_e$ by a limit of a sum of the terms indexed by
geometrical objects called clusters where the number of terms can be estimated 
and for each term, a large deviation technique with an appropriate scaling can be
applied. The main steps of our proof can be summarized as follows. 
 
\noindent
 {\it1. Proof of the inequality} 
\begin{align}
\log r^*_e \leq \limsup_{a\to +\infty}  &\limsup_{b\to +\infty} \limsup_{T\to +\infty}   \nonumber\\
&\sup_{\substack{z : \; az \in E, \\ N <
  |z| \leq b N}} \frac{1}{aT} \log \P_{az} \bigl(|Z_a(T)|\leq bN, \; |Z_a(s)| \geq N, \; \forall \;
0\leq s \leq T\bigr) \label{e3-1}
\end{align}
for $N>0$ large enough. $\P_{az}(\cdot)$ denotes here the conditional probability given that $Z_a(0) =
z$ (or equivalently $X(0)=az$). If the order of the limits in
$a$ and $T$ could be  reversed,  a large deviation upper bound would
give  directly  a  good  estimate  of  the right  hand  side  of  this
inequality.  Unfortunately, such an inversion of the limits 
  seems  very difficult  to  prove :  for this one should be able to perform large deviation estimates
simultaneously for all $T$ large enough with the same large deviation parameter $a$ large
enough.   An alternative approach  consists in sub-dividing the interval $[0,T]$ on smaller disjoint
intervals $[t_{i-1},t_i]$ is such a way that for every of these intervals, the desired
large deviation estimates can be performed with the same parameter $a$ large enough. For
this we use the following preliminary results.
 
\medskip
 
\noindent
{\it 2. Some preliminary results and constructions.} 
Let $M_N(t,x,y)$ denote the infimum of the rate function $I_{[0,t]}(\phi)$ over all
$\phi:[0,t]\to{\cal E}$ such that $\phi(0)=x$, $\phi(t)=y$ and $$\sup_{0\leq s\leq
  t}|\phi(s)| \geq N.$$ Using the upper large deviation bound we show that for given $\eps > 0$ and
$N>0$ there exists a finite set $V(\eps,N) \subset \{(x,y) : |x|\leq 2N, \; |y|\leq 2N\}$
  and there are strictly positive real numbers $\delta(x,y)$ and $a(x,y)$ such that  
\[
\{(x,y) : |x|\leq 2N, \; |y|\leq 2N\} \subset \bigcup_{(x,y)\in V(\eps,N)}
B(x,\delta(x,y)) \times B(y, \delta(x,y))
\]
and for any $(x,y)\in
  V(\eps,N)$, for all $a\geq a(x,y)$,  and for any $z\in \frac{1}{a} E$ satisfying the  
inequality $|z-x|\leq \delta(x,y)$, the following inequalities hold
\begin{align*}
\log\P_{az}\left(\sup_{s\in [0,1]} |Z_a(s)|
\geq N  \; \mbox{ and } \; Z_a(1)\in
B\bigl(y,\delta(x,y)\bigr)\right)
&\leq - a M_N(1,x,y) + a\eps
\end{align*}
when $M_N(1,x,y) <+\infty$, and
\begin{equation*}
\log\P_{az}\left( \sup_{s\in [0,1]} |Z_a(s)|
\geq N \; \mbox{ and } \; Z_a(1) \in B\bigl(y,\delta(x,y)\bigr)
\right)
\leq - a I^*_0 N/\eps
\end{equation*}
when $M_N(1,x,y) =+\infty$. Here and throughout, $B(x,\delta)$ denotes the open ball
centered at $x$ and having radius $\delta$.

\medskip
 
\noindent
{\it 3. Change of scale}. Using the above estimates and the identity $Z_{a}(s) =  t
    Z_{at}(s/t)$ we conclude that for any $(x,y)\in
  V(\eps,N)$, for all $\sigma > 0$, $a\geq \sigma a(x,y)$ and $t\geq 1/\sigma$,  and for any $z\in
    \frac{1}{a} E\cap t B(x,\delta(x,y))$, the following inequalities hold
\begin{multline}\label{e4-1}
\log\P_{az}\left(\sup_{s\in [0,t]} |Z_a(s)|
\geq Nt  \; \mbox{ and } \; Z_a(t)\in
t B\bigl(y,\delta(x,y)\bigr)\right) \\
\leq - a M_N(t,tx,ty) + at\eps
\end{multline}
when $M_N(t,tx,ty) <+\infty$, and
\begin{equation}\label{e5-1}
\log\P_{az}\left( \sup_{s\in [0,t]} |Z_a(s)|
\geq N t\; \mbox{ and } \; Z_a(t) \in t B\bigl(y,\delta(x,y)\bigr)
\right)
\leq - a t I^*_0 N/\eps
\end{equation}
otherwise. These inequalities hold simultaneously for all $(x,y)\in
  V(\eps,N)$ and $t\geq 1/2 $ when $a > 2 \max_{(x,y)\in V(\eps,N)} a(x,y)$. Such a change of
  scale is a key point of our proof. 
 
\medskip
 
\noindent
{\it 4. Cluster expansion}. For given $\eps>0$, $N>0$ and $T$, and for every function
    $\phi\in D([0,T],{\cal E})$ with $
 |\phi(t)| \geq N$ for all $0\leq t\leq T$ and such that $|\phi(0)| \leq 2TN$ and
    $|\phi(T)|\leq 2TN$, 
we define a partition $0=t_0 < t_1 < \ldots < t_n = T$ and a sequence $\bigl((x_i,y_i)\in  
    V(\eps,N), \; i=1,\ldots,n\bigr)$, as follows. 
\begin{itemize}
\item[--] If $\sup_{0\leq t\leq T} |\phi(t)| > NT$ we let $n=1$ and we choose
  $(x_1,y_1)\in V(\eps,N)$ such that $\phi(0) \in T B(x_1,\delta(x_1,y_1))$ and $\phi(T)
  \in T B(y_1,\delta(x_1,y_1))$.
\item[--] Otherwise, we divide the interval $[0,T]$ in two intervals $[0,T/2]$ and
  $[T/2,T]$ and we restart our construction for the restriction of 
  $\phi$ on each of them. 
\end{itemize}
This algorithm terminates because $|\phi(t)| \geq N$ for all $0\leq t\leq T$. The
resulting sequence $\Gamma(\phi) =
\bigl( (t_1,x_1,y_1),\ldots, (t_n,x_n,y_n)\bigr)$ is called 
$(T,\eps, N)$- cluster corresponding to $\phi$ (see section~\ref{sec6} for a more careful
definition of $(T,\eps, N)$- cluster). In statistical physics, the notion of cluster is
usually associated with a connected graph. In our context, the cluster $\Gamma(\phi)$ is
connected in the following sense :
\[
\phi(t_i) \in  (t_i-t_{i-1})
B(y_i, \delta(x_i,y_i)) \cap (t_{i+1}-t_i) B(x_{i+1}, \delta(x_{i+1},y_{i+1})) \not= \emptyset
\]
 for every $1\leq i\leq n$. For each cluster $\Gamma =
\bigl( (t_1,x_1,y_1),\ldots, (t_n,x_n,y_n)\bigr)$ we consider the quantity 
\[
\chi_a(\Gamma) = \sup_{z}
P_{az}\bigl( \Gamma(Z_a) = \Gamma\bigr)
\]
where the supremum is taken over all $z\in \frac{1}{a}E \cap t_1 B(x_1,
\delta(x_1,y_1))$. Using inequality (\ref{e3-1}) we obtain 
\[
\log r^*_e \leq \limsup_{a\to +\infty} \limsup_{b\to +\infty} \limsup_{T\to +\infty}
\frac{1}{aT}  \log\left(\sum_\Gamma \chi_{a}(\Gamma) \right) 
\]
where the summation is taken over all clusters $\Gamma =$ 
$((t_1,x_1,y_1),\ldots, (t_n,x_n,y_n))$ for
which the sets $t_1 B(x_1,
\delta(x_1,y_1) \cap B(0,bN)$ and $(t_n-t_{n-1}) B(y_n,
\delta(x_n,y_n)\cap B(0,bN)$ are 
non-empty.  
 
\medskip
 
\noindent
{\it 4. Cluster estimates}. We 
show that for given $\eps >0$, $N>0$ and $T\geq 1$, there are at 
most $(2|V(\eps,N)|)^{2T}$  clusters. For every cluster $\Gamma =$ 
$((t_1,x_1,y_1),\ldots, (t_n,x_n,y_n))$ we obtain a good estimate of the
quantity $\chi_{a}(\Gamma)$ by using the inequalities  (\ref{e4-1}),
(\ref{e5-1}) with $x=x_i$, $y=y_i$ and $t=t_i-t_{i-1}$ for every $i=1,\ldots,n$ (see
Lemma~\ref{lem6-6} below) and we deduce from them the  desired inequality $\log r^*_e \leq - I^*_0$.
 
\bigskip

Our paper is organized as follows.  Section~\ref{secII} is devoted to general properties 
of the quantities $r^*$ and $r^*_e$.  In section~\ref{sec4},  
 the definition and some useful properties of sample path large
deviations are recalled and  different
representations of the quantities $I^*$ and $I^*_0$ are derived. Section~\ref{sec5}  is devoted to
the proof of Theorem~\ref{th1}. Theorem~\ref{th2} is proved
in section~\ref{sec6}. In section~\ref{sec7} we apply our results to calculate the
quantities $r^*$ and $r^*_e$ for two dimensional Jackson networks.

To simplify the notations, we consider continuous time Markov processes. For discrete time Markov
processes our results can be extended in a straightforward way.
 
 \section{Some general properties}\label{secII}
In this section, we recall general properties of the quantities $r^*$ and
$r^*_e$.   
 
\subsection{Spectral radius}\label{sec2}
For $r>0$ and for a real bounded measurable function $\varphi$ on $E$, the resolvent
function is defined by   
\[
R_r\varphi(x) = \int_0^\infty r^{- t} P^t\varphi(x) \, dt
\]
where $P^t$ denotes Markov semi-group corresponding to the process $(X(t))$~:
$$P^t\varphi(x) =\E_x(\varphi(X(t)).$$
Recall that by definition,  $r^*$ 
 is the infimum of all  those $r>0$ for which the function $R_r\1_V$ is $m$-integrable on
 compact subsets of $E$ for every compact set $V\subset E$.   
 
It is clear moreover that $r^*$ 
 is the infimum of all  those $r>0$ for which the function $R_r\varphi$ is $m$-integrable on
 compact subsets of $E$ for every non-negative continuous function $\varphi$ having a
 compact support. 
 
The following property of the quantity 
$r^*$ immediately
follows from the definition. 
\begin{prop}\label{pr0-2}  
\begin{equation}\label{e1-2}
\log r^* \leq  \sup_{W,V} \; \limsup_{t\to\infty} \; \sup_{x\in
  W}  \, \frac{1}{t} \log  P^t\1_V(x)
\end{equation}
where the supremum $\sup_{W,V}$ is taken over all compact subsets $W,V\subset E$. 
Moreover, if for a compact set $V\subset E$, there exists a compact set $V'\subset E$ and there
are  real numbers $t>0$ and $\eps >0$  such that
$P^s\1_{V'} \geq \eps\1_V$ for all $0<s<t$, then for any compact subset $W\subset E$, 
\begin{equation}\label{e3-2}
\log r^* \geq   \limsup_{t\to\infty} \; \inf_{x\in
  W} \, \frac{1}{t} \log   P^t\1_V(x). 
\end{equation}
\end{prop} 
\begin{proof} Indeed, Relation (\ref{e1-2}) holds because  by Fubini's theorem 
\[
m(\1_W R_r\1_V) ~=~ \int_0^\infty r^{-t} m(\1_W P^t\1_V(x) \, dt  ~\leq~ m(W) \int_0^\infty r^{-t} 
\sup_{x\in
  W} P^t\1_V(x) \, dt.  
\]
Suppose moreover that for a compact set $V\subset E$ there
are  real numbers $t>0$ and $\eps >0$ and a compact set $V'\subset E$ such that
$P^s\1_{V'}(x) \geq \eps$ for all $0<s<t$ and for all $x\in V$. 
Then  for any  
increasing sequence of real positive numbers $t_n$ with $\inf_n(t_{n+1} -t_n) \geq t$, the
following inequalities hold 
\begin{align*}
m(\1_WR_r\1_{V'}) &\geq~ \sum_n \int_{t_n}^{t_n+t} r^{-s} \,
m(\1_WP^{s}\1_{V'}) \, ds \\&\geq~ \sum_n \int_{0}^{t} r^{-t_n-s} \, 
m(\1_WP^{t_n} \1_V P^s\1_{V'}) \, ds \\ &\geq~ \eps \sum_n r^{-t_n}  \, m( \1_WP^{t_n} \1_V)
\int_{0}^{t} r^{-s} \, ds.  
\end{align*}
The last inequality shows that for $x\in W$,  $\, m(\1_WR_r\1_{V'}) = +\infty$ whenever  
\[
0< r < \limsup_n \, \frac{1}{t_n} \log  m(\1_W P^{t_n}\1_V)   
\]
and consequently, Inequality (\ref{e3-2}) is verified.
\end{proof}
 
The following proposition shows that for a large class of Markov processes, the quantity
$r^*$ can be represented in terms of $r$-superharmonic functions. 
 
\begin{defi}
 A measurable function $f : E\to \R_+$ is called  $r$-superharmonic with $r>0$ if the
inequality $P^t f \leq r^t f$ holds for all $t\in\R_+$. We say that a function $f$ is
locally $m$-integrable on $E$ if it is $m$-integrable on the compact subsets of $E$. 
\end{defi} 
\begin{prop}\label{pr1-2} Suppose that  there exists a non-negative continuous function
  $\varphi_0$ on $E$ having a compact support such that for every $t>0$ the function
  $P^t\varphi_0$ is continuous on $E$ and for every $x\in E$ there exists
  $t>0$ such that $P^t\varphi_0(x) > 0$.  Then  $r^*$  is
  the infimum of all those
  $r>0$ for which there exists a non-negative locally $m$-integrable $r$-superharmonic function
  $f$ with $\inf_{x\in W} f(x) > 0$ for every compact subset $W\subset E$. 
\end{prop} 
\begin{proof} Indeed, for $r>r^*$, the function
  $R_r\varphi_0$ is non-negative locally $m$-integrable and $r$-superharmonic. Moreover, the
  sample paths of the Markov process $(X(t))$ being right continuous, the mapping $t\to
  P^t\varphi_0(x)$ is right continuous and hence, under the hypotheses
  of our proposition,  $R_r\varphi_0(x) > 0$ for every $x\in E$.  Finally, by Fatou's
  lemma, the function $R_r\varphi_0$  is lower semi-continuous on $E$ and consequently,
  for any compact set $W\subset E$, 
\[
\inf_{x\in W} R_r\varphi_0(x) > 0. 
\]
Conversely, suppose that for $r>0$ there exists a non-negative locally $m$-integrable
  $r$-superharmonic function $f$ 
such that $\inf_{x\in W} f(x) > 0$ for every compact subset $W\subset E$. Then for every
compact set $W\subset E$, and for any $t>0$ and $x\in W$, 
\[
P^t\1_W (x)  \inf_{y\in W} f(y)   ~\leq~ P^t f(x) \leq  r^t f(x).  
\]
> From this it follows that for any $r'>r$, the function $R_{r'}\1_W$ is locally $m$-integrable and
consequently,  $r^* \leq r$.   Proposition~\ref{pr1-2} is therefore
verified.  
\end{proof}
 
\subsection{Essential spectral radius}\label{sec3}
Recall that by definition, the quantity $r^*_e$ is the infimum over all
those $r>0$ for which there is a compact set
$K\subset E$, such that the the truncated resolvent function 
\begin{equation}\label{e1-3}
R_{K,r}\1_W(x) = \int_0^\infty r^{-t} \P_x(X(t)\in W, \, \tau(K) > t)  \,dt 
\end{equation}
is $m$-integrable on the compact subsets of $E\setminus K$ for every compact set $W\subset
E\setminus K$. 
 
\begin{prop}\label{pr0-3} 
\begin{equation}\label{e2-3}
\log r^*_e ~\leq~  \inf_K \; \sup_{W,V\subset E\setminus K} \; \limsup_{t\to\infty} \; \sup_{x\in
  W} \, \frac{1}{t} \log   \P_x(X(t)\in V, \, \tau(K) > t)     
\end{equation}  
where the infimum is taken over all compact subsets $K\subset E$ and the supremum
is taken over all compact subsets $W,V\subset E\setminus K$.  
Suppose moreover that for any compact subset $K\subset E$ there are compact sets $V_K,V'_K\subset
E\setminus K$ and there is a real number $t>0$ such that 
\[
\inf_{0<s\leq t} \; \inf_{x\in V_K} \P_x(X(s)\in V'_K, \; \tau(K) > t) ~>~0. 
\]
Then 
\begin{equation}\label{e3-3}
\log r^*_e ~\geq~ \inf_K \; \sup_{W\subset E\setminus K} \; \limsup_{t\to\infty} \;
\inf_{x\in W}  \, 
\frac{1}{t} \log   \P_x( X(t) \in V_K, \, \tau(K) >t)   
\end{equation}
where the infimum $\inf_K$ is taken over all compact subsets $K\subset E$ and the supremum
is taken over all compact subsets $W\subset E\setminus K$.  
\end{prop} 
\begin{proof} Indeed, let $r^*_K$ be the infimum over all those $r>0$ for which the
  function  (\ref{e1-3}) is $m$-integrable on the compact subsets of $E\setminus K$ for every
  compact set $W\subset E\setminus K$. Using the same
  arguments as for the proof of Proposition~\ref{pr0-2} we obtain that for every compact
  set $K\subset E$,   
\[
\log r^*_K ~\leq~  \sup_{W,V\subset E\setminus K} \; \limsup_{t\to\infty}
\; \sup_{x\in W}  \, \frac{1}{t} \log  \P_x( X(t) \in V, \, \tau(K) >t) 
\]
and moreover, for compact subsets $W,V\subset E\setminus K$, 
\[
\log r^*_K ~\geq~ \limsup_{t\to\infty} \; \inf_{x\in W}  \,
\frac{1}{t} \log  \P_x( X(t) \in V_K, \, \tau(K) >t)  
\]
if there are a compact set $V'\subset E\setminus K$ and a real numbers $t>0$  such that
\[
\inf_{0<s\leq t} \; \inf_{x\in V} \P_x(X(s)\in V', \; \tau(K) > s) ~>~ 0. 
\] 
Using relation $r^*_e = \inf_K r^*_K$, this proves Proposition~\ref{pr0-3}.  
\end{proof}

The next proposition describes the quantity $r^*_e$ in terms of Lyapunov
functions which are superharmonic outside of compact sets. 
 
\begin{defi} A measurable function $f : E\to \R_+$ is called  $r$-superharmonic outside of
a compact set $K\subset E$ with $r>0$ if the
inequality 
\[
\E_x(f(X(t)) ; \; \tau(K) > t)  \leq r^t f(x)
\]
holds for all $x\in E\setminus K$ and for all $t\in\R_+$.
\end{defi}
\begin{prop}\label{pr1-3} Suppose that  there exists an increasing sequence of open
  relatively compact sets $U_n\subset E$ such that $\cup_n U_n = E$ and let for every
  $n\in\N$ there exist  a non-negative continuous function
  $\varphi_n$ on $E$ having a compact support in $E\setminus\ol{U}_n$ such that for every
  $t>0$ the function $
x\to E_x(\varphi_n(X(t)) ; \; \tau(\ol{U}_n) > t) $ 
is continuous on $E\setminus\ol{U}_n$ and for every $x\in E\setminus\ol{U}_n$ there exists
  $t>0$ such that $$E_x(\varphi_n(X(t)) ; \; \tau(\ol{U}_n) > t) > 0.$$  Then  $r^*_e$  is
  the infimum of all those $r>0$ for which there exists a compact set 
  $K\subset E$ and a non-negative locally $m$-integrable $r$-superharmonic  outside of $K$
  function $f$ with $\inf_{x\in W} f(x) > 0$ for every compact subset $W\subset E\setminus K$.  
\end{prop} 
\begin{proof} Given a compact subset $K\subset E$, let  $r^*_K$ be the infimum over all
  those $r>0$ for which the function  (\ref{e1-3}) is $m$-integrable on the compact subsets of
  $E\setminus K$ for every compact set $W\subset E\setminus K$, and let $\rho_K$ be the
  infimum of all those $r>0$ for which there exists  a non-negative locally $m$-integrable
  $r$-superharmonic  outside of $K$  function $f$ with $\inf_{x\in W} f(x) > 0$ for every
  compact subset $W\subset E\setminus K$.  The same arguments as in the proof of
  Proposition~\ref{pr1-2} show that for $K_n = \ol{U}_n$, $r^*_{K_n} = \rho_{K_n}$ for every
  $n\in\N$. The quantities $r^*_K$ and $\rho_K$ being decreasing with respect to $K$, this
  implies that 
\[
r^*_e = \inf_n r^*_{K_n} = \inf_n \rho^*_{K_n} = \inf_K \rho_K
\]
where the last infimum is taken over all compact subsets $K\subset E$.
Proposition~\ref{pr1-3} is therefore proved.  
\end{proof}

\begin{defi} Let $\sigma^*_e$ denote 
the infimum over all those $r>0$ for which there is a
  compact set $K\subset E$ such that the function  
\begin{equation}\label{e4-3}
x \to \int_0^{+\infty} r^{-t} \, \P_x (\tau(K) > t)  \,dt 
\end{equation}
is $m$-integrable on the compact subsets of $E\setminus K$.
\end{defi}
 
The following proposition represents the number $\sigma^*_e$ in terms of
Lyapunov functions. 
 
\begin{prop}\label{pr2-3}  $\sigma^*_e$  is the infimum of all those $r>0$ for which
  there exists a compact set 
  $K\subset E$ and a non-negative locally $m$-integrable $r$-superharmonic  outside of $K$
  function $f$ with $\inf_{x\in E\setminus K} f(x) > 0$. 
\end{prop}
\begin{proof}  Suppose that for $\hat{r}>0$, there is a non-negative $\hat{r}$-superharmonic outside
  of $K$ function $f$  which is $m$-integrable on the compact subsets of $E\setminus K$ and such
  that $\inf_{x\in E\setminus K} f(x) > 0$. Then for any $t>0$, and for any
  $x\in E\setminus K$, 
\[
  \P_x (\tau(K) > t)  \times \inf_{y\in E\setminus K} f(y) ~\leq~ \E_x( f(X(t) ; \; \tau(K) > t) ~\leq~ 
\hat{r}^t f(x). 
\]
For any $r>\hat{r}$, the function (\ref{e4-3}) 
is therefore $m$-integrable on the compact subsets of $E\setminus K$ and consequently, $\hat{r}\geq
\sigma^*_e$. The function
  (\ref{e4-3}) being $r$-superharmonic outside of $K$, this
proves Proposition~\ref{pr2-3}. 
\end{proof}
 
The next proposition shows that for a large class of recurrent Markov processes, the
quantities $r^*_e$ and $
\sigma^*_e$ are equal.  
 
\begin{prop}\label{pr3-3}
Suppose that the hypotheses of Proposition~\ref{pr1-3} are satisfied and let 
for every $n>0$,   and  $\tau(\ol{U}_{n+1}) < \tau(\ol{U}_n) < +\infty$, $\P_x$-almost
surely for every $x\in E\setminus \ol{U}_{n+1}$. Then $r^*_e =  \sigma^*_e$. 
\end{prop}
\begin{proof}
If for $r>0$ there is a compact set $K\subset E$ such that 
  the function (\ref{e4-3}) is $m$-integrable on the compact subsets of $E\setminus K$, 
 then for every compact set $W\subset E$, the function (\ref{e1-3}) 
is $m$-integrable on the compact subsets of $E\setminus K$ and therefore, $\sigma^*_e \geq
r^*_e$. 
 
Let us prove that $\sigma^*_e \leq r^*_e$. Because of Proposition~\ref{pr2-3},  it is sufficient
to show that for any $r>r^*_e$ there is $n\in\N$ and a non-negative $r$-superharmonic
outside of $K_n=\ol{U}_n$ function 
$f$ which is $m$-integrable on the compact subsets of $E\setminus K_n$ and such that $\inf_{x\in
  E\setminus K_n} f(x) > 0$. It is sufficient moreover to consider the case when $r^*_e <
r<1$ because for $r\geq 1$, the function $f\equiv 1$ is $r$-superharmonic. 
 
Given a compact subset $K\subset E$, let  $r^*_K$ be the infimum over all
  those $r>0$ for which the function  (\ref{e1-3}) is $m$-integrable on the compact subsets of
  $E\setminus K$ for every compact set $W\subset E\setminus K$. The quantity $r^*_K$ being
  decreasing with respect to $K$,
\[
r^*_e = \inf_n r^*_{K_n}.
\]
Given  $r^*_e < r<1$ let us choose $k\in\N$ such that $r >
r^*_{K_k}$ and let $n>k$ be such that $K_k\subset U_{n-1}$.  Under the hypotheses of
Proposition~\ref{pr1-3} the function 
\[
f(x) = \int_0^{+\infty} r^{-t} \E_x(\varphi_k(X(t)) ; \; \tau(K_k) > t)  \, dt 
\]
is $m$-integrable on the compact subsets of $E\setminus K_k$ and $r$-superharmonic
outside of $K_k$. Moreover the same arguments as for the proof of Proposition~\ref{pr1-2}
show that for any compact subset $W\subset E\setminus K_k$, there is $\eps(W) >0$ such
that 
\[
\inf_{x\in W} f(x) ~\geq~ \eps(W).
\]
Consider  a compact set $W_n=K_{n}\setminus U_{n-1}$. Using
strong Markov property, we obtain  
\begin{align*}
\inf_{x\in E\setminus K_n} f(x) &~\geq~  \inf_{x\in E\setminus K_n} \; \int_0^{+\infty}
r^{-t} \E_x\bigl(\varphi_k(X(t)) ; \;  \tau(K_k) > t \geq \tau(W_n)\bigr)  \, dt  \\ &~=~
\inf_{x\in E\setminus K_n} \; \E_x\Bigl(r^{-\tau(W_n) }
f (X(\tau(W_n)) ) ; \; \tau(K_k) > \tau(W_n) \Bigr)  \\ &~\geq~ \eps(W_n) \, \inf_{x\in
  E\setminus K_n} \; \P_x\Bigl(\tau(K_k) > \tau(W_n) \Bigr) ~=~ \eps(W_n)   
\end{align*}
The last equality holds here because 
under the hypotheses of our proposition, for all $x\in E\setminus K_n$, $\P_x$-almost
surely  $\tau(K_n) = \tau(W_n) < 
\tau(K_{n-1}) < \tau(K_k)$. The function $f$ being $m$-integrable on the compact subsets of
$E\setminus K_n$ and $r$-superharmonic outside of $K_n$ the last relation shows that 
 $\sigma^*_e \leq r$. Letting finally $r\to r^*_e$  we conclude that $\sigma^*_e \leq
r^*_e$. The equality $\sigma^*_e = r^*_e$ is 
therefore proved. \end{proof}

 \section{Sample path large deviations}\label{sec4}
\subsection{Definitions and general properties}

Let $D([0,T],\R^d)$ denote the set of all right continuous with left
limits functions from $[0,T]$ to $\R^d$ endowed with Skorohod metric
$d_S(\cdot,\cdot)$. Recall that Skorohod metric $d_S(\phi,\psi)$ is defined as the 
infimum of those positive $\eps$ for which there exists a strictly
increasing continuous mapping $\lambda$ from $[0,T]$ onto itself 
satisfying inequalities 
\[
\sup_{t>s} \left| \log \frac{\lambda(t)-\lambda(s)}{t-s} \right|
\leq \eps \quad \mbox{ and } \quad
\sup_t |\phi(t)-\psi(\lambda(t))| \leq \eps
\]
where the metric is induced by the Euclidean norm $|\cdot|$ on $\R^d$.
 
The space $D([0,T],\R^d)$ endowed with Skorohod
metric is complete. A sequence
$\phi_n\in D([0,T],\R^d)$ converges to a limit $\phi\in D([0,T],\R^d)$ in the
Skorohod metric if and only if there exist strictly increasing
continuous mappings $\lambda_n: [0,T]\to [0,T]$ such that
$\lambda_n(0)=0$, $\lambda_n(T)= T$, $\lambda_n(t)\to t$ as
$n\to\infty$ uniformly in $t\in [0,T]$ and $\phi_n\circ\lambda_n(t)
\to \phi(t)$ as $n\to\infty$ uniformly in $t\in[0,T]$. When $\phi$ is
continuous on $[0,T]$, Skorohod convergence $\phi_n\to\phi$ implies
uniform convergence.  For non-continuous $\phi$, Skorohod
convergence $\phi_n\to\phi$ implies $\phi_n(0)\to \phi(0)$,
$\phi_n(T)\to\phi(T)$ and $\phi_n(t)\to\phi(t)$ for continuity points
$t\in [0,T]$ of $\phi$. For more details about Skorohod metric, we refer the reader
to Billingsley~\cite{Billingsley}.
 
We consider a Markov process $(X(t))$ on
$E\subset 
\R^d$.  The trajectories of the Markov process $X(t)$ are assumed to be 
almost surely right continuous and to have  left limits. 
 
\smallskip
\noindent
\begin{defi} 
1) A mapping $I_{[0,T]}:~D([0,T],\R^d)\to
[0,+\infty]$ is a {\it good rate function} on $D([0,T],\R^d)$ if for
any $c>0$ and for any compact set $V\subset \R^d$, the set
\[
\{ \varphi \in D([0,T],\R^d): ~\phi(0)\in V \; \mbox{
and } \; I_{[0,T]}(\varphi) \leq c \}
\]
is compact in $D([0,T],\R^d)$.  According to this definition, a good
rate function is lower semi-continuous. 
 
2) The family of scaled Markov processes $$(Z_a(t), \; t\in[0,T])\stackrel{\text{def.}}{=}
(X(at)/a, \; t\in[0,T])$$ is said to 
satisfy {\it sample path large deviation principle} in $D([0,T], \R^d)$ with a rate function
$I_{[0,T]}$ if for any $x\in\R^d$ 
\begin{equation}\label{e-lb-4}
\lim_{\eps\to 0} \;\liminf_{a\to\infty} \; \inf_{y: |y-x|<\eps} \frac{1}{a}\log\P_{ay}\left( Z_a(\cdot)\in 
{\cal
O}\right) \geq -\inf_{\phi\in{\cal O}:\phi(0)=x} I_{[0,T]}(\phi), 
\end{equation}
for every open set ${\cal
O}\subset D([0,T],\R^d)$,
and
\begin{equation}\label{e-ub-4}
\lim_{\eps\to 0} \;\limsup_{a\to\infty} \; \sup_{y
  : |y-x|<\eps} \frac{1}{a}\log\P_{ay}\left( Z_a(\cdot)\in
F\right) \leq -\inf_{\phi\in F:\phi(0)=x} I_{[0,T]}(\phi).
\end{equation}
 for every closed set $F\subset   D([0,T],\R^d)$. 
\end{defi} 
$\P_{ay}$ denotes here and throughout a conditional probability given that $Z_a(0) = y$ 
(or equivalently $X(0)=ay$), the infimum at the left hand side of the inequality
(\ref{e-lb-4}) and the supremum at the  left  hand side of the inequality
(\ref{e-ub-4}) are taken over all $y\in \frac{1}{a} E$ satisfying the inequality $|x-y| <
\eps$. 
 
\smallskip
 
We refer to sample path large deviation
principle as SPLD principle. Inequalities (\ref{e-lb-4}) and (\ref{e-ub-4}) are
referred as lower and upper SPLD bounds respectively.
 
The following statement is a consequence of contraction principle (see Dembo and
Zeitouni~\cite{D-Z}) and the identity $Z_{a}(t) =  T Z_{aT}(t/T)$.
 
\begin{prop}\label{pr1-4}
Suppose that the family of scaled processes
$Z_a(t)$ satisfies SPLD principle in $D([0,1],\R^d)$ with a
good rate functions $I_{[0,1]}$.
Then for any
$T>0$, this family of processes satisfies SPLD principle in
$D([0,T],\R^d)$ with the good rate functions
\begin{equation}\label{e1-4}
I_{[0,T]}(\phi) = T I_{[0,1]}(G_T\phi) 
\end{equation}
where $G_T$ denotes  a mapping from $D([0,T],\R^d)$ to $D([0,1],\R^d)$
defined by
\[
G_T\phi(t) = \phi(Tt)/T, \quad t\in [0,1].
\]
\end{prop}
When the SPLD principle holds, the corresponding rate function is unique and hence, under
the hypotheses~(A), relation (\ref{e1-4}) is  satisfied. A good rate function
satisfies moreover the following properties.

\medskip
 
\begin{lemma}\label{lem1-4}
For any closed set $F\subset D([0,T],\R^d)$,
the mapping 
\begin{equation}\label{e2-4}
(x,y) \to \inf_{\substack{\phi\in F:~ \phi(0)=x,\\ \phi(T)=y}}
I_{[0,T]}(\phi)
\end{equation}
is lower semi-continuous on $\R^{2d}$. 
\end{lemma}
\begin{proof}
Indeed, for any $c>0$  and for any compact set $V\subset\R^d$, 
\begin{equation}\label{e3-4}
\{ (x,y)\in V\times\R^{d} :  \inf_{\substack{\phi\in F:~ \phi(0)=x,\\
\phi(T)=y}} I_{[0,T]}(\phi) \leq c \} = \\ \bigcap_n \xi(K_n) 
\end{equation}
where $K_n = \{ \phi\in F : \phi(0)\in V, \; I_{[0,T]}(\phi) \leq c + 1/n \}$ and the
mapping $\xi : D([0,T],\R^d)\to\R^2$ is defined by  
$\xi(\phi) = (\phi(0),\phi(T))$.
The sets $K_n$
being compact and the mapping $\xi$ being continuous, the sets
$\xi(K_n)$ are compact.  This proves that the set \eqref{e3-4} is also compact and
consequently, the mapping \eqref{e2-4} is lower semi-continuous. 
\end{proof}

\subsection{Different representations of the quantities $I^*$  and
$I^*_0$.}\label{sec4p} Throughout
this section, $I_{[0,T]}$ denotes a good rate
function on $D([0,T],\R^d)$ satisfying Assumption (A) and Relation (\ref{e1-4}). Recall
that ${\cal E}$ denotes the set of all possible limits $\lim_{a\to\infty} x_a/a$ with
$x_a\in E$. According to Assumption (A), 
\begin{itemize}
\item[--] the set ${\cal E}\setminus\{0\}$ is convex and
non-empty; 
\item[--] the rate function $I_{[0,T]}$ has an integral form : there is a local rate function
$L:{\cal E}\times\R^d\to\R_+$ such that 
\[
I_{[0,T]}(\phi) = \int_0^T L(\phi(t),\dot\phi(t))  \, dt 
\]
if the function $\phi : [0,1]\to {\cal E}$ is absolutely continuous, and $I_{[0,1]}(\phi)
= +\infty$ otherwise.
\item[--] there is a convex function  $l_2$
  on $\R^d$  that is finite in a neighborhood of zero and  such that  
  $L(x,v) \leq l_2(v)$ for all $x\in{\cal E}$ and for
  all $v\in\R^d$. 
\end{itemize}
Here and throughout, for  $x,y\in{\cal E}$ and $t>0$,  $I(t,x,y)$ denotes the infimum of the rate function
$I_{[0,t]}(\phi)$ over  all continuous functions $\phi : [0,t]\to{\cal E}$ with
$\phi(0)=x$ and $\phi(t)=y$, $I_0(t,x,y)$ denotes the infimum
of the rate function $I_{[0,t]}(\phi)$ over all continuous functions 
$\phi : [0,t]\to{\cal E}$ such that
$\phi(0)=x$, $\phi(t)=y$ and $\phi(s)\not=0$ for all $s\in]0,t[$, and
$\hat{I}(t,x,y)$ is the infimum of the rate function $I_{[0,t]}(\phi)$ over
all those $\phi : [0,t]\to{\cal E}$ for which $\phi(0)=x$, $\phi(t)=y$ and
the set $\{s\in [0,t] : \phi(s)=0\}$ has Lebesgue measure zero.
By definition, 
\[
 I^* =  I(1,0,0) \quad \text{ and  } \quad  I^*_0 = \hat{I}(1,0,0).
\]
 
\noindent
{\bf Remark 1}. Because of Relation (\ref{e1-4}),  for any $x,y\in{\cal E}$ and 
$t>0$, 
\begin{equation}\label{e1-4p}
I(t, tx, ty) = t I(1, x, y), \quad \quad \hat{I}(t, tx, ty) = t \hat{I}(1, x, y), 
\end{equation}
and 
\begin{equation}\label{e2-4p}
I_0(t, tx, ty) = t I_0(1, x, y).  
\end{equation}
\noindent
{\bf Remark 2}. The integral representation of the rate function implies that for any
$\phi\in D([0,T],\R^d)$ and  $0<t<T$, 
\begin{equation}\label{e3-4p}
I_{[0,T]}(\phi) = I_{[0,t]}(\phi) + I_{[t,T]}(\phi)
\end{equation}
where $I_{[t,T]}(\phi) = I_{[0,T-t]}(\phi(t+\cdot))$. From this it follows that for any $x,y,z\in{\cal E}$,
\[
I(T,x,y) \leq I(t,x,z) + I(T-t,z,y) \; \text{ and } \; \hat{I}(T,x,y)
\leq \hat{I}(t,x,z) + \hat{I}(T-t,z,y).
\]
Moreover, for $z\not= 0$, 
\[
I_0(T,x,y) \leq I_0(t,x,z) + I_0(T-t,z,y).
\]
 
\medskip 
 
\noindent
{\bf Remark 3}. Since the function $l_2$ is convex and finite in a neighborhood of zero, there
are $c > 0$ and $C>0$ such that $l_2(v)\leq C$ for all $v\in\R^d$ with $|v|\leq c$. Hence,
for $\phi(t) = x + t(y-x)$ with $x,y\in{\cal E}$ such that $|x-y| \leq c$,   
\[
I(1,x,y) ~\leq~ \hat{I}(1,x,y) ~\leq~ I_0(1,x,y) ~\leq~ I_{[0,1]}(\phi) ~\leq~ l_2(y-x) \leq C.
\] 
Using Relations (\ref{e1-4p}) and (\ref{e2-4p}) this implies that 
\begin{equation}\label{e4-4p}
I(t,x,y) ~\leq~ \hat{I}(t,x,y) ~\leq~ I_0(t,x,y) ~\leq~ Ct 
\end{equation}
for all $x,y\in{\cal E}$ and for all $t > |y-x|/c$. 
\medskip
 
Using Relations (\ref{e2-4p})  and  (\ref{e3-4p}) we obtain the following statement. 
\begin{lemma}\label{lem2-4} 
For any continuous function $\phi : [0,T]\to {\cal E}$ with $\phi(0) = \phi(T) = 0$,
\begin{equation}\label{e5-4p}
I_{[0,T]}(\phi) \geq I_{[0,T]}(\ol{0})
\end{equation}
and 
\begin{equation}\label{e6-4p}
I_{[0,T]}(\phi) \geq    I_0(1,0,0) \, \mes\{t\in[0,T]:\phi(t)\not=0\} 
\end{equation}
where $\mes\{t\in[0,T]:\phi(t)\not=0\}$ denotes Lebesgue measure of the set
$\{t\in[0,T]:\phi(t)\not=0\}$. 
\end{lemma}
\begin{proof} Indeed, given a continuous function $\phi : [0,T]\to {\cal E}$ with $\phi(0)
  = \phi(T) = 0$,  consider a sequence of functions $\phi_n : [0,T]\to {\cal E}$ defined by $ 
\phi_n(t+{k}/{n}) = \phi(nt)/n$ for all $
t\in [0,{T}/{n}]$ and for all $k=0,\ldots,n$.
It is clear that $\phi_n\to \ol{0}$ as $n\to\infty$ uniformly on
$[0,T]$. The rate function $I_{[0,T]}(\cdot)$ being lower
semi-continuous, this implies that
\[
\liminf_{n\to\infty} I_{[0,T]}(\phi_n) \geq I_{[0,T]}(\ol{0})
\]
because  uniform convergence implies Skorohod
convergence. Relations
(\ref{e1-4}) and (\ref{e3-4p}) show  
that $I_{[0,T]}(\phi_n) = n I_{[0,T/n]}(\phi_n) = I_{[0,T]}(\phi)$
and hence, the last inequality  proves relation (\ref{e5-4p}).
 
Furthermore, for such a function $\phi$, the set
  $\{t\in]0,T[:\phi(t)\not=0\}$ is a union of countable family of open disjoint intervals
$]t_i,t'_i[ \subset [0,T]$, $i\in\N$. Using relation (\ref{e3-4p})  we obtain
\[
I_{[0,T]}(\phi) \geq \sum_i I_{[t_i,t'_i]}(\phi) \geq \sum_i
I_0(t'_i-t_i,0,0) = I_0(1,0,0) \sum_i (t'_i-t_i)
\]
where the second relation holds because  for every $i$,
$\phi(t_i)=\phi(t'_i)=0$ and $\phi(t)\not= 0$ for all
$t\in]t_i,t'_i[$, and the third relation follows from 
relation (\ref{e1-4}). The last relation proves inequality (\ref{e6-4p}). 
\end{proof}

The following proposition gives several equivalent representations of $I^*$ and $I^*_0$. 
 
\begin{prop}\label{prop1-4} 1) For any $x,y\in{\cal E}$ and $T>0$, 
\begin{equation}\label{e7-4p}
I^* ~=~ I_{[0,1]}(\ol{0}) ~=~ \limsup_{t\to +\infty} I(t,x,y)/t  ~=~  \inf_{\phi(0)=\phi(T)} 
I_{[0,T]}(\phi)/T
\end{equation}
where the infimum is taken over all continuous functions $\phi : [0,T]\to {\cal E}$ with
$\phi(0)=\phi(T)$,  and  
\begin{align}
I^*_0 ~=~ I_0(1,0,0) &=~ \limsup_{t\to +\infty} I_0(t,x,y)/t  ~=~ \limsup_{t\to +\infty}
\hat{I}(t,x,y)/t \nonumber \\ &=~ \inf_{\substack{\phi(0)=\phi(T),  \\
\phi(t)\not= 0, \; \forall \; 0<t<T} }  I_{[0,T]}(\phi)/T \label{e8-4p}
\end{align}
where the infimum is taken over all continuous functions $\phi : [0,T]\to {\cal E}$ with
$\phi(0)=\phi(T)$ such that $\phi(t)\not=0$ for all $0<t<T$. 
\end{prop}
 
\begin{proof}
The equalities $I^* = I_{[0,1]}(\ol{0})$ and $I^*_0 = I_0(1,0,0)$ hold because of Lemma~\ref{lem1-4}.
Furthermore, Relation (\ref{e4-4p}) and the inequality 
$$I(t_1+t_2+t, x,y) \leq I(t_1,x,x') + I(t,x,y) + I(t_2,y',y)$$ 
show that the limits $\; \liminf_{t\to\infty}  I(t,x,y)/t \, $ and $\;
\limsup_{t\to\infty}  I(t,x,y)/t\, $ do not depend on $x,y\in {\cal E}$. Since by
(\ref{e1-4p}), $I(t,0,0) = t I(1,0,0) = I^*$, we conclude that for any
$x,y\in {\cal E}$
\[ I(t,x,y)/t  \to  I(1,0,0) \quad \text{  as } \; t\to\infty. 
\]
Moreover, Relation (\ref{e3-4p}) shows that $I(nT,x,x) \leq nI(T,x,x)$ for all $x\in{\cal
  E}$, $T>0$ and $n\in\N$ and consequently,  
\[
\frac{1}{T} I(T,x,x) \geq \lim_{n\to +\infty} \frac{1}{nT} I(nT,x,x) = I^*.
\]
The last relation combined with the inequality 
\[
I^* \; = \;  I(T,0,0)/T \; \geq \;
 \inf_{x\in{\cal E}} I(T,x,x)/T \; = \;
 \inf_{\phi(0)=\phi(T)} I_{[0,T]}(\phi)/T
\]
proves that 
\[
I^*  ~=~  \inf_{x\in{\cal E}} I(T,x,x)/T ~=~  \inf_{\phi(0)=\phi(T)} I_{[0,T]}(\phi)/T.
\]
Relation (\ref{e7-4p}) is therefore proved. The same arguments show that for all $x\in{\cal
  E}$ and $T>0$
\[
I^*_0 ~=~ \limsup_{t\to +\infty} \hat{I}(t,x,y)/t  ~=~  \inf_{\substack{\phi(0)=\phi(T),
    \\ \mes\{ t : \phi(t)= 0\} = 0} }  I_{[0,T]}(\phi)/T 
\]
where the infimum is taken over all continuous functions $\phi : [0,T]\to{\cal E}$ with
$\phi(0)=\phi(T)$ for which the set $\{t : \phi(t) =0\}$ has Lebesgue measure zero, that  
the limits $ l = \liminf_{t\to +\infty} I_0(t,x,y)/t$, $ l_1 = \liminf_{t\to +\infty} I_0(t,0,y)/t$
and $ l_2 = \liminf_{t\to +\infty} I_0(t,x,0)/t$ do not depend on
$x,y\in{\cal E}\setminus\{0\}$ and that for every $T>0$,
\begin{align*}
I_0(1,0,0) \; &= \; I_0(T,0,0)/T \; = \lim_{t\to +\infty} I_0(t,0,0)/t  \\
&\leq  \; l_i \; \leq \; l \; \leq \; \inf_{x\not=0} I_0(T,x,x)/T \; = \;
\inf_{\substack{\phi(0)=\phi(T)\not=0,\\
\phi(t)\not= 0, \forall t\in]0,T[} } I_{[0,T]}(\phi)/T, \quad i=1,2.
\end{align*}
Hence, to complete the proof of Relation (\ref{e8-4p}) it is sufficient to show that the
infimum at the right hand side of the above relation does not exceed
$ I_0(T,0,0)/T$, or equivalently that  
\begin{equation}\label{e8-4pp}
I_{[0,T]}(\phi) \geq \inf_{x\not=0} I_0(T,x,x).
\end{equation}
for any continuous function $\phi :[0,T]\to{\cal E}$ with $\phi(0)=\phi(T)=0$ and
$\phi(t)\not=0$ for all $t\in]0,T[$. For such a function $\phi$, for every $0 < \delta < T/2$ 
 there is a function
$\psi_\delta : [0,t_\delta]\to {\cal E}$ with  $t_\delta = |\phi(T - \delta) -
 \phi(\delta)|/c$,  $\psi_\delta(0) = \phi(T-\delta)$ and 
 $\psi_\delta(t_\delta) = \phi(\delta)$ such that
$\psi_\delta(t)\not=0$ for all $t\in[0,t_\delta]$ and
$
I_{[0,t_\delta]}(\psi_\delta) \leq  Ct_\delta$ (see Remark 3).
Define the function $\phi_\delta$ by
setting $\phi_\delta(t) = \phi(t+\delta)$ for $t\in[0,T-2\delta]$ and
$\phi_\delta(t) = \psi_\delta(t-T+2\delta)$ for
$t\in[T-2\delta,T-2\delta + t_\delta]$. Then, using relation (\ref{e3-4p}), we get
\[
I_{[0,T - 2\delta + t_\delta]}(\phi_\delta) = I_{[0,t_\delta]}(\psi_\delta) +
I_{[\delta,T -\delta]}(\phi) \leq Ct_\delta + I_{[0,T]}(\phi)
\]
and using relation (\ref{e2-4p}) it follows that
\begin{align*}
I_{[0,T - 2\delta + t_\delta]}(\phi_\delta) &\geq
I_0(T - 2\delta + t_\delta,\phi(\delta),\phi(\delta)) \\ &\geq \inf_{x\not=0}
I_0(T - 2\delta + t_\delta,x,x) = \frac{1}{T}(T - 2\delta + t_\delta)
\inf_{x\not=0} I_0(T,x,x).
\end{align*}
This proves that
\[
I_{[0,T]}(\phi) \geq \frac{1}{T}(T - 2\delta + t_\delta) \inf_{x\not=0} I_0(T,x,x) -
Ct_\delta.
\]
Letting at the last inequality $\delta\to 0$ we obtain Relation (\ref{e8-4pp}).
 
\end{proof}

 \section{Proof of Theorem~\ref{th1}.}\label{sec5}
Recall that 
$I(t,x,y)$  denotes
the infimum of the rate functions $I_{[0,t]}(\phi)$  over all continuous functions
$\phi : [0,t]\to{\cal E}$, with $\phi(0)=x$ and $\phi(t)=y$.
By Proposition~\ref{prop1-4},
\[
I^* ~=~ I_{[0,1]}(\ol{0}) ~=~ \limsup_{t\to +\infty} I(t,x,y)/t  ~=~  \inf_{\phi(0)=\phi(T)} 
I_{[0,T]}(\phi)/T
\]
where the infimum is taken over all continuous functions $\phi : [0,T]\to {\cal E}$ with
$\phi(0)=\phi(T)$. To complete the proof of Theorem~\ref{th1} it is sufficient therefore
to show that $\log r^* = -I^*$. 
 
To prove the upper bound 
\begin{equation}\label{e1-5}\log
r^* \leq -I^*. 
\end{equation}
we use Proposition~\ref{pr0-2} and the SPLD upper
bound. The SPLD upper  bound implies that for all compact sets $V,W\in{\cal K}$,  
\begin{align*}
\limsup_{t\to\infty} \; \sup_{x\in
  W}  \, \frac{1}{t} \log  P^t\1_V(x) &\leq~ \lim_{\delta\to0}
\lim_{\eps\to0} \limsup_{a\to\infty} 
\frac{1}{a} \log  \sup_{x : |x| < \eps}  \P_{ax}\bigl( |Z_a(1)| \leq
\delta \bigr) \\
&\leq~ - \lim_{\delta\to0}  \; \inf_{\phi(0)=0, 
\; |\phi(1)| \leq \delta} I_{[0,1]}(\phi) 
\end{align*}
where the infimum is taken over all continuous functions $\phi :[0,1]\to {\cal E}$ with
$\phi(0)=0$ and $|\phi(1)|\leq \delta$. The right hand side of the last inequality equals 
\[
- \lim_{\delta\to0}  \; \inf_{y : \; |y| \leq \delta} I(1,0,y). 
\]
The mapping $(x,y) \to I(1,x,y)$ being lower semi-continuous (see Lemma~\ref{lem1-4}),
using Proposition~\ref{pr0-2} we conclude that 
\[
\log r^* ~\leq~  \sup_{W,V} \; \limsup_{t\to\infty} \; \sup_{x\in
  W}  \, \frac{1}{t} \log  P^t\1_V(x) ~\leq~ -I(1,0,0) ~=~ -I^*. 
\]
Inequality (\ref{e1-5}) is therefore proved. 
 
\bigskip
 
To prove the lower bound $\log
r^* \geq -I^*$ we use Relation (\ref{e3-2}) of Proposition~\ref{pr0-2} and the SPLD lower
bound. The hypotheses of Proposition~\ref{pr0-2} are satisfied for every compact set
$V\subset E$ with $t=a$ and $V' = V'_{a\delta} = \{x\in E: |x|\leq a \delta\}$ for 
$\delta >0$ and for $a>0$  large enough because   
\[
\inf_{0<s<a} \; \inf_{y\in V} \; P^s\1_{V'}(y) \geq   \inf_{y:
  |y|<\eps} \ \P_{ay}\bigl( \|Z_a(s)\|_\infty < \delta  \bigr)
\]
if $V\subset\{x : |x| < a\eps\}$, and because by the SPLD lower bound,
\begin{align*}
 \lim_{\eps\to 0} \; \liminf_{a\to\infty} \; \inf_{y: |y|<\eps}
\frac{1}{a}\log\P_{ay}\bigl( \|Z_a(s)\|_\infty < \delta \bigr) &\geq~ -\inf_{\phi
  :\phi(0) = 0, \|\phi\|_\infty < \delta } I_{[0,T]}(\phi) \\ &\geq~ - I_{[0,1]}(\ol{0}) ~\geq - l_2(0)
~> -\infty. 
\end{align*}
Using Relation (\ref{e3-2})  with $W= V=\{x\in E : |x|\leq a\eps\}$ and letting $a\to
+\infty$, we obtain  
\[
\log r^* \geq 
 \liminf_{a\to \infty}  \liminf_{n\to \infty} 
\frac{1}{an} \inf_{|x|\leq \eps} 
  \log \P_{ax}\bigl( |Z_a(n)| < \eps \bigr). 
\]
Hence, by Markov property,  
\begin{align}
\log r^* &\geq~  \liminf_{a\to \infty} \; \liminf_{n\to+\infty}  
 \;  \inf_{|x| < \eps} \; \frac{1}{an} \log
 \P_{ax}\bigl( |Z_a(k)| < \eps, \; \forall \; 1\leq k
\leq n \bigr) \nonumber \\
&\geq~  \liminf_{a\to \infty} \; 
\inf_{|x| < \eps} \; \frac{1}{a} \log \P_{ax}\bigl(
 |Z_a(1)| < \eps \bigr).  \label{e2-5}
\end{align}
Furthermore, SPLD lower bound proves that for any $\sigma >0$ and for any
$x\in\R^d$, there exist $\delta(x)>0$ and $a(x)>0$ such
that for all $a\geq a(x)$ and for all $x'\in \frac{1}{a} E$ satisfying inequality $|x'-x| < \delta(x)$, 
\begin{align*}
\frac{1}{a} \log
\P_{ax'}\left( |Z_a(1)| < \eps \right) &\geq
-\inf_{\phi :~\phi(0)=x, |\phi(T)|<\eps} I_{[0,T]}(\phi) - \sigma \\ &\geq
-\inf_{y : |y| < \eps} I(1,x,y) - \sigma  ~\geq~ - I(1,x,0) - \sigma.  
\end{align*}
Recall  that $E\subset{\cal E}$ and consequently, $\frac{1}{a} E\subset {\cal E}$.  
The set $\{x\in{\cal E} : |x|\leq\eps\}$ being compact, there are
$x_1,\ldots,x_n\in\{x\in {\cal E} : |x| \leq \eps\}$ such that
\[
\{x\in {\cal E} : |x|\leq\eps\}\subset \bigcup_{i=1}^n \{x':~ |x'-x_i| <\delta(x_i)\}.
\]
For $a\geq \max\{a(x_1),\ldots,a(x_n)\}$, we obtain therefore 
\begin{align*}
\inf_{x\in\frac{1}{a} E : |x| < \eps} \; \frac{1}{a} \log
\P_{ax}\left( |Z_a(1)| < \eps  \right)  &\geq - \max_i I(1,x_i,0) - 
\sigma  \geq 
- \sup_{x : |x| \leq \eps}  I(1,x,0) - 
\sigma.   
\end{align*}
The last relation combined
with Inequality (\ref{e2-5}) proves that  
\begin{equation}\label{e3-5}
\log r^* \geq - \sup_{x : |x| \leq \eps}  I(1,x,0).  
\end{equation}
Finally, Relations~(\ref{e1-4p}) and (\ref{e3-4p}) show that for any $0<t<1$, 
\[
I(1,x,0) \leq I(t,x,0) +
I(1-t,0,0) = I(t,x,0) + (1-t)I(1,0,0) \leq  I(t,x,0) + I^*\]
and consequently, for $t=\eps/c$, using Inequality (\ref{e4-4p}) it follows that  
\begin{multline*}
\log r^* \geq  - \sup_{x : |x| \leq \eps}  I(1,x,0) \geq  - \sup_{x : |x| \leq \eps}
I(\eps/c , x, 0) - I^*  = - C\eps/c - I^*.   
\end{multline*}
Letting at the last inequality $\eps \to 0$, we conclude  $
\log r^* \geq - I^*$.  
The inequality $\log r^* \geq - I^*$ combined with  (\ref{e1-5}) proves that $\log r^* = - I^*$.

 \section{Proof of Theorem~\ref{th2}}\label{sec6}
 
Recall that  for $x,y\in\R^d$ and $t>0$,  $\hat{I}(t,x,y)$ denotes 
the infimum of $I_{[0,t]}(\phi)$ over all continuous functions
$\phi :[0,t]\to{\cal E}$ with $\phi(0)=x$, $\phi(t)=y$ for which the set $\{s\in[0,t]~: \phi(s)=0 \}$
has Lebesgue measure zero. $I_0(t,x,y)$ denotes
the infimum of $I_{[0,t]}(\phi)$ over all continuous functions
$\phi :[0,t]\to{\cal E}$ such that $\phi(0)=x$, $\phi(t)=y$  and
$\phi(s)\not=0$ for all $0<s<t$.  By Proposition~\ref{prop1-4}, 
\begin{align}
I^*_0 ~=~ I_0(1,0,0) &=~ \limsup_{t\to +\infty} I_0(t,x,y)/t  ~=~ \limsup_{t\to +\infty}
\hat{I}(t,x,y)/t  \label{e0-6} \\ &=~ \inf_{\substack{\phi(0)=\phi(T),  \\
\phi(t)\not= 0, \; \forall \; 0<t<T} }  I_{[0,T]}(\phi)/T \nonumber 
\end{align}
where the infimum is taken over all continuous functions $\phi : [0,T]\to {\cal E}$ with
$\phi(0)=\phi(T)$ such that $\phi(t)\not=0$ for all $0<t<T$. 
 
To prove Theorem~\ref{th2} it is sufficient therefore to show that 
$\log r^*_e = - I^*_e$.

\subsection{Lower bound $\log r^*_e \geq -I^*_0$.} The proof of  this inequality uses 
Relation (\ref{e3-3}) of Proposition~\ref{pr0-3} and SPLD lower bound. We begin our proof
with the following lemma.
 
\begin{lemma}\label{lem1-6p} Under the hypotheses~(A), for any compact subset $K\subset
  E$, for any  $x\in{\cal E}\setminus\{0\}$, and for any $\delta > \eps >0$ 
such that $|x| >\delta$, the inequality 
\begin{equation}\label{e1-6}
\inf_{0<s\leq t} \; \inf_{x\in V_K} \P_x(X(s)\in V'_K, \; \tau(K) > s) ~>~ 0. 
\end{equation}
holds with $V_K =
\{y\in E : |ax-y| \leq a\eps\}$ and $V'_K = \{y\in E : |ax-y| \leq a\delta\}$ for all
$a>0$ large enough. 
\end{lemma}
\begin{proof} Let $x\in{\cal E}\setminus\{0\}$, $\sigma > 0$ and $\delta > \eps >0$ be
such that $|x| >\delta + \sigma$. For a compact subset $K\subset E$ are is $a_K >0$
such that for $a>a_K$, $K\subset \{ y : |y| < a \sigma \}$. For $V_K =
\{y\in E : |ax-y| \leq a\eps\}$ and $V'_K = \{y\in E : |ax-y| \leq a\delta\}$,  
\begin{multline}\label{e2-6}
\inf_{0<s\leq a} \; \inf_{x\in V_K} \P_x(X(s)\in V'_K, \; \tau(K) > s) ~\geq \\\inf_{y
  : |y-x|<\eps} \P_{ay}\left( \sup_{0\leq s \leq 1} |Z_a(s) - x| < \delta\right)  
\end{multline}
for all $a>a_K$. Moreover because of SPLD lower bound and Assumption~$(a_3)$, 
\begin{align*}
\lim_{\eps\to 0} \;\liminf_{a\to\infty} \; \inf_{y: |y-x|<\eps}
\frac{1}{a}\log\P_{ay}\left( \sup_{0\leq s \leq 1} |Z_a(s) - x| < \delta\right) &\geq -
I_{[0,1]}(\phi) \\ &\geq~ -l_2(0) ~>~ -\infty 
\end{align*}
where $\phi(t)\equiv x$. The last inequality shows that for  any 
 $\eps > 0$ small enough, there is $a_\eps > a_K$ such that  for every $a>a_\eps$, the right hand side of
 Inequality (\ref{e2-6}) is strictly positive and hence,
Inequality (\ref{e1-6}) holds. 
\end{proof}
 
Because of Lemma~\ref{lem1-6p} and 
Proposition~\ref{pr0-3},  Relation  (\ref{e3-3})  holds for any compact set
 $K\subset E$ with $V_K = \{y\in E : |ax-y| \leq a\eps\}$ for any $0<\eps<|x|$ 
 and for all $a>0$ large enough.  Using this relation together with Markov property we obtain 
\begin{align*}
\log \, &r^*_e ~\geq~ \inf_K \; \limsup_{t\to\infty} \; \inf_{x\in V_K}  \, 
\frac{1}{t} \log   \P_x( X(t) \in V_K, \, \tau(K) >t)  \\  
\geq&~ \liminf_{a\to +\infty} 
\liminf_{n\to \infty} \inf_{y : |y-x| \leq \eps} \frac{1}{anT} \log \P_{ay}
\left(|Z_a(nT)-x| < \eps, \inf_{s\in[0,nT]} |Z_a(s)| > \sigma\right) \\ 
\geq&~ \liminf_{a\to +\infty} 
\lim_{n\to \infty} \inf_{y : |y-x| \leq \eps} \frac{1}{anT} \log \P_{ay}
\left(\max_{k\leq n} |Z_a(kT)-x| < \eps, \inf_{s\in[0,nT]} |Z_a(s)| > \sigma\right) \\ 
\geq&~ \liminf_{a\to \infty} \inf_{y : 
  |y-x| \leq \eps}\frac{1}{aT} \log \P_{ay} \left(|Z_a(T)-x| < \eps, \; \inf_{s\in[0,T]} |Z_a(s)| > \sigma 
\right)
\end{align*} 
for any $T>0$.  The last inequality combined with the SPLD lower bound, the 
same arguments as in the proof of Theorem~\ref{th1} (see the proof of inequality
(\ref{e2-5})) and Relation (\ref{e3-4p}) show that for any  $T > 1$,  
\begin{align}
\log r^*_{e} &\geq~ - \sup_{y\in{\cal E} : |x-y| \leq \eps} \; 
I_\sigma(T,y,x)/T \nonumber \\ &\geq~ - \sup_{y\in{\cal E}  : |x-y| \leq \eps} \; I_\sigma(1,y,x)/T ~-~ 
I_\sigma(T-1,x,x)/T  \label{e3-6}  
\end{align}
where $I_\sigma(t,x,y)$ denotes the infimum of the rate function
$I_{[0,t]}(\phi)$ over all those $\phi\in D([0,t],{\cal E})$
for which  $\phi(0)=y$, $\phi(t)=x$ and $|\phi(s)| > \sigma$ for all
$s\in[0,t]$. Moreover, letting $\phi_y(t) = y + (x-y) t$ we get  
\[|\phi_y(t)| \geq |x| - |\phi_y(t)-x|
\geq |x| - \eps > \sigma \]
for all $t\in[0,1]$ and consequently, by Assumption~(A), 
\[
I_\sigma(1, y,x) ~\leq~ I_{[0,1]}(\phi_y) ~=~ \int_0^1 L(\phi_y(t),
\dot\phi_y(t))  \, dt ~\leq~ l_2(y-x).   
\]
Using this relation for the right hand side of Inequality~(\ref{e3-6}) and letting
 $\sigma\to 0$ we obtain 
\[
\log r^*_{e} \geq - \sup_{v : |v| \leq \eps}  l_2(v) /T - 
I_0(T-1,x,x)/T. \]
The function $l_2$ being finite in a neighborhood of zero, this implies that 
\[
\log r^*_{e} \geq - \limsup_{T\to\infty} I_0(T-1,x,x)/T.
\] 
The last inequality combined with Relation~(\ref{e0-6}) proves that $\log r^*_{e} \geq -
I^*_0$.

 \medskip 
 
\subsection{Upper bound $\log r^*_e \leq - I^*_0$.} 
To prove this relation we use Proposition~\ref{pr0-3} and the SPLD upper bound. Because of
Proposition~\ref{pr0-3}, for any $N>0$, 
\begin{multline} \label{e1-7}
\log r^*_e \leq \limsup_{a\to +\infty} \, \lim_{b\to +\infty} \, 
\limsup_{T\to +\infty}   \\   
\sup_{N <
  |z| \leq b N}\frac{1}{aT} \log \P_{az} \bigl(|Z_a(T)|\leq bN, \; |Z_a(s)| \geq N, \; \forall \;
0\leq s \leq T\bigr). 
\end{multline}
To estimate the right hand side of this inequality, we use cluster expansion method and
SPLD upper bound. Before to introduce the notion of clusters, we consider the following
preliminary results. 
 
 For $x,y\in{\cal E}$ and $t>0$, let  $M_N(t,x,y)$ denote the infimum of the rate
function $I_{[0,t]}(\phi)$ aver all functions 
$\phi \in D([0,t],{\cal E})$ such that   $\phi(0)=x$, $\phi(t)=y$ and 
\[
\sup_{0\leq s \leq t}|\phi(s)| \geq N t.
\] 
 
\begin{lemma}\label{lem2-6} For  
$x,y\in{\cal E}$ and $N>0$, 
\begin{multline}\label{e3-7}
\lim_{\delta\to 0}\limsup_{a\to\infty} \sup_{z\in 
  B(x,\delta)} \;
\frac{1}{a}\log\P_{az}\Bigl(\sup_{s\in [0,1]}
|Z_a(s)|\geq N  \mbox{ and }  Z_a(1) \in B(y,\delta) \Bigr) \\ 
 \leq \; - M_N(1,x,y).
\end{multline}
\end{lemma}
\begin{proof} Using SPLD upper bound it follows
that for any $x,y\in {\cal E}$ and for any $\sigma>0$,
\begin{align*}
&\lim_{\delta\to 0}\limsup_{a\to\infty} \sup_{z\in B(x,\delta)} \; \frac{1}{a}\log\P_{az}\Bigl(\sup_{s\in 
[0,1]} 
|Z_a(s)|\geq N \; \mbox{ and } \; Z_a(1)\in B(y,\delta)\Bigr) \\
&\; \leq \lim_{\delta\to 0}\limsup_{a\to\infty} \sup_{z\in
  B(x,\delta)}  \; \frac{1}{a}\log\P_{az}\Bigl(\sup_{s\in  
[0,1]} |Z_a(s)|\geq N \; \mbox{ and } \;
Z_a(1)\in\ol{B}(y,\sigma)\Bigr)
\\
&\; \leq  - \inf_{y'\in\ol{B}(y,\sigma)}
  M_N(1,x,y').
\end{align*}
By Lemma~\ref{lem1-4}, the
mapping $(x,y)\to M_N(1,x,y)$ is lower semi-continuous and hence, 
letting at the last relation $\sigma\to 0$ we get inequality~(\ref{e3-7}).
\end{proof}
 
\medskip
 
 By Lemma~\ref{lem2-6}, for any $\eps>0$, $N>0$ and  
 $x,y\in{\cal E}$, there exist $\delta(x,y)>0$ and $a(x,y)>0$  such
that for any
$a>a(x,y)$ and for any $z\in \frac{1}{a}E$ satisfying the  
inequality $|z-x|\leq \delta(x,y)$, the following inequalities hold. 
\begin{align*}
\log\P_{az}\left(\sup_{s\in [0,1]} |Z_a(s)|
\geq N  \; \mbox{ and } \; Z_a(1)\in
B\bigl(y,\delta(x,y)\bigr)\right)
&\leq - a M_N(1,x,y) + a\eps
\end{align*}
when $M_N(1,x,y) <+\infty$, and
\begin{equation*}
\log\P_{az}\left( \sup_{s\in [0,1]} |Z_a(s)|
\geq N \; \mbox{ and } \; Z_a(1) \in B\bigl(y,\delta(x,y)\bigr)
\right)
\leq - a I^*_0 N/\eps
\end{equation*}
when $M_N(1,x,y) =+\infty$. Moreover,  the real number
$a\geq a(x,y)$ at the above inequalities can be replaced by $at$ with $a\geq 2
a(x,y)$ and $t\geq 1/2$. Using relation $Z_{at}(s) = Z_a(ts)/t$ this implies  the following statement.
 
\begin{lemma}\label{lem3-6} For any $\eps >0$ and $N>0$, there are strictly
positive functions $\delta(\cdot,\cdot)$ and $a(\cdot,\cdot)$ on ${\cal E}\times {\cal E}$
such that for any $x,y\in{\cal E}$,
$t>1/2$, $a>2 a(x,y)$ and $z\in \frac{1}{a}E$ 
  satisfying the inequality $|z-tx|\leq t\delta(x,y)$, the following
  inequalities hold
\begin{multline}\label{e4-7}
\log\P_{az}\left(\sup_{s\in [0,t]} |Z_a(s)|
\geq N t  \; \mbox{ and } \; Z_a(t)\in
t B\bigl(y, \delta(x,y)\bigr)\right)  \\ \leq - a t M_N(1,x,y) +
at\eps
\end{multline}
when $M_N(1,x,y) <+\infty$, and
\begin{equation}\label{e5-7}
\log\P_{az}\left(\sup_{s\in [0,t]} |Z_a(s)|
\geq N t \; \mbox{ and } \; Z_a(t) \in t B\bigl(y, \delta(x,y)\bigr)
\right)
\leq - a t I^*_0 N/\eps
\end{equation}
otherwise.
\end{lemma}
 
We are ready now to introduce the notion of cluster. 
For a given $\eps > 0$, we choose
$N\geq 1$ large enough so that 
\begin{equation}\label{e2-7}
I_{[0,t]}(\phi) >  I^*_0/\eps
\end{equation}
for any function 
$\phi \in D([0,t],{\cal E})$ for which there are $0 < s_1 < s_2 \leq 1$ such that
$
|\phi(s_1) - \phi(s_2) | \geq N
$. 
Such a choice is possible because under the hypotheses~(A), for any absolutely continuous
function $\phi : [0,1] \to {\cal E}$ for which there are $0 < s_1 < s_2 \leq 1$ such that
$
|\phi(s_1) - \phi(s_2) | \geq N,
$ 
\begin{align*}
I_{[0,1]}(\phi) &\geq~ \int_{s_1}^{s_2} L(\phi(s),\dot\phi(s)) ds ~\geq~ \int_{s_1}^{s_2}
l_1(\dot\phi(s)) ds \\ &\geq~ (s_2-s_1) l_1\bigl((\phi(s_2) - \phi(s_1))/(s_2-s_1)\bigr)
~\geq~ N \inf_{v : |v| \geq N} l_1(v)/|v| ~\to~ +\infty  
\end{align*}
when $N\to\infty$. 
\smallskip
 
Given $\eps>0$ and $N\geq 1$, let
$\delta(\cdot,\cdot)$ and $a(\cdot,\cdot)$ be the positive functions 
on ${\cal E}\times{\cal E}$ satisfying Lemma~\ref{lem3-6}.  
Without any restriction of 
generality we will suppose that for any
$x,y\in{\cal E}$, 
\begin{equation}\label{e6-7}
\delta(x,y)\leq \varepsilon. 
\end{equation}  
The set $\ol{B}(0,2N)\times\ol{B}(0,2N)$ being compact there exists a
finite subset $V(N,\eps)\subset \ol{B}(0,2N)\times\ol{B}(0,2N)$ such
that
\[
\ol{B}(0,2N)\times\ol{B}(0,2N) \subset \bigcup_{(x,y)\in V(N,\eps)}
B(x,\delta(x,y))\times B(y,\delta(x,y)).  
\]
For a function $\phi\in D([T_1,T_2],{\cal E})$ such that  $
 |\phi(t)| \geq N$ for all $T_1\leq t\leq T_2$ and $|\phi(T_1)| \leq 2(T_2-T_1)N$,  
    $|\phi(T_2)|\leq 2(T_2-T_1)N$, we define a partition $$T_1=t_0 < t_1 < \ldots <
    t_n = T_2$$ and a sequence $\bigl((x_i,y_i)\in V_{(\eps,N)}, \; i=1,\ldots,n\bigr)$ by
    induction : 
\begin{itemize}
\item[--] We  choose
  $(x_1,y_1)\in V(\eps,N)$ such that $\phi(T_1) \in (T_2-T_1) B(x_1,\delta(x_1,y_1))$ and $\phi(T_2)
  \in (T_2-T_1) B(y_1,\delta(x_1,y_1))$, and we let $n=1$ if 
\[
\sup_{T_1\leq t\leq T_2} |\phi(t)| > N(T_2-T_1).
\]
\item[--] Otherwise, we divide the interval $[T_1,T_2]$ in two intervals $[T_1,(T_2+T_1)/2]$ and
  $[(T_2-T_1)/2,T_2]$ and we restart our construction for the restriction of 
  $\phi$ on each of them. 
\end{itemize}
This algorithm terminates because $|\phi(t)|\geq N$ for all $t\in[T_1,T_2]$. The resulting
    sequence $$\Gamma(\phi) =
\bigl( (t_1,x_1,y_1),\ldots, (t_n,x_n,y_n)\bigr)$$ is called 
{\it $(\eps, N)$- cluster corresponding to the function $\phi$}. Remark that for such a
    sequence, 
\[
\phi(t_i) \in (t_i-t_{i-1})B\bigl(y_i,\delta(x_i,y_i)\bigr) \cap (t_{i+1}-t_i)
B\bigl(x_{i+1},\delta(x_{i+1},y_{i+1})\bigr)\not= \emptyset
\] 
 for every $i=1,\ldots,n-1$ and $(t_i-t_{i-1}) \geq 1/2$ for all $i=1,\ldots,n$. Moreover,
 for $T_1=0$ and $T_2=T>1$ and for a natural number $k$ such that $2^{k-1} < T
 \leq 2^k$,  by construction, $t_i \,2^k/T \in \N$ for all $i=0,\ldots,n$. 
 
\smallskip
\noindent 
\begin{defi}
 A sequence $\Gamma =
\bigl( (t_1,x_1,y_1),\ldots, (t_n,x_n,y_n)\bigr)$ is called {\it $(T,\eps,N)$- cluster} if
\begin{itemize}
\item[--]  $0=t_0 < t_1 < \ldots < t_n =T$ and $(t_i-t_{i-1}) \geq 1/2$ for all  $i=1,\ldots,n$; 
\item[--]  $t_i \,2^k/T \in \N$ for all $i=0,\ldots,n$;   
\item[--]  the set $(t_i-t_{i-1})B\bigl(y_i,\delta(x_i,y_i)\bigr) \cap (t_{i+1}-t_i)
B\bigl(x_{i+1},\delta(x_{i+1},y_{i+1})\bigr)$ is non-empty for every $i=1,\ldots,n-1$. 
\end{itemize} 
\end{defi}
For a given $(T,\eps, N)$- cluster $\Gamma =
\bigl( (t_1,x_1,y_1),\ldots, (t_n,x_n,y_n)\bigr)$, we denote by $D(\Gamma)$ the
set of all functions $\phi
:[0,T]\to{\cal E}$ such that $\Gamma(\phi) = \Gamma$. The quantity $\chi_a(\Gamma)$ is defined by  
\[
\chi_a(\Gamma) = \sup_{x}
P_x\bigl( Z_a\in D(\Gamma)\bigr)
\]
where the supremum is taken over all $x\in E$ such that
$|x - a t_1 x_1| \leq a t_1 \delta(x_1,y_1)$.
 
Remark that by construction,
\begin{equation}\label{e7-7}
\sup_{N <
  |z| \leq b N} P_{az}\bigl(|Z_a(T)|\leq bN, \; |Z_a(s)| \geq N, \; \forall \;
0\leq s \leq T\bigr) \, \leq \, \sum_\Gamma \chi_{a}(\Gamma),   
\end{equation} 
where the sum is taken over all $(T,\eps,N)$-clusters $\Gamma =$ 
$((t_1,x_1,y_1),\ldots, (t_n,x_n,y_n))$ for
which the sets $t_1 B(x_1,
\delta(x_1,y_1) \cap B(0,bN)$ and $(t_n-t_{n-1}) B(y_n,
\delta(x_n,y_n)\cap B(0,bN)$ are 
non-empty.
 
To estimate the right hand side of the inequality (\ref{e7-7}) (and hence also the right
hand side of the inequality (\ref{e1-7})) we estimate the number of $(T,\eps,N)$-clusters
and for each $(T,\eps,N)$-cluster $\Gamma$ we estimate the quantity
$\chi_{a}(\Gamma)$. This is a subject of the following lemmas. 
 
\begin{lemma}\label{lem4-6} For any $T\geq 1$, there are at most 
$(2|V(N,\eps)|)^{2T}$ $\; (T,\eps,N)$-clusters, where $|V(N,\eps)|$ denotes the number of points of
the set $V(N,\eps)$.
\end{lemma}
\begin{proof} Indeed, let $2^{k-1} < T\leq
  2^k$. Then there are at most $2^{2^k}\leq 2^{2T}$ possible
choices for a sequence of natural numbers $0=t_0 \,2^k/T <t_1 \,2^k/T <\ldots <t_n \,2^k/T =2^k$
and for any $1\leq n\leq 2^k$ there are at most $|V(N,\eps)|^{n}\leq
|V(N,\eps)|^{2^{k}}\leq |V(N,\eps)|^{2T}$ possible choices for a sequence $(x_1,y_1),\ldots,
(x_n,y_n)\in V(N,\eps)$.
\end{proof} 
 
To estimate the quantities $\chi_{a}(\Gamma)$ we use Assumption~(A) and
Lemma~\ref{lem3-6}. Recall that under the hypotheses~(A), the rate function $I_{[0,T]}$
has an integral form~: there is a local rate function $L:{\cal E}\times\R^d\to\R_+$ such that 
\[
I_{[0,T]}(\phi) = \int_0^T L(\phi(t),\dot\phi(t))  \, dt 
\]
if the function $\phi : [0,1]\to {\cal E}$ is absolutely continuous, and $I_{[0,1]}(\phi)
= +\infty$ otherwise. Moreover, there is a  finite in a
neighborhood of zero convex function $l_2 : \R^d\to\R_+$ such that $L(x,v) \leq l_2(v)$ for
all $x\in{\cal E}$ and for all $v\in\R^d$.  
The function $l_2$ being convex and finite in a neighborhood of zero, there are real
numbers $C>0$ and $c>0$ such that 
\begin{equation}\label{e8-7}
\sup_{v :|v|\leq c} l_2(v) \leq C.
\end{equation}
 
\begin{lemma}\label{lem6-6} For any $T > b>1$, for any $a\geq
  2 \max\{a(x,y) \; | \; (x,y)\in V(N,\eps)\}$ and for any $(T,\eps,N)$-cluster $\Gamma =
  \bigl((t_1,x_1,y_1),\ldots,(t_n,x_n,y_n)\bigr)$ for
which the sets $t_1 B(x_1,
\delta(x_1,y_1) \cap B(0,bN)$ and $(t_n-t_{n-1}) B(y_n,
\delta(x_n,y_n)\cap B(0,bN)$ are 
non-empty,  the following inequality holds
\begin{equation}\label{e9-7}
\frac{1}{a}\log\chi_a(\Gamma) \leq  - T(1-2\eps) I^*_0 + 4T\eps C/c + T\eps + 2bN C/c
\end{equation}
 \end{lemma}
\begin{proof}  Using 
  Lemma~\ref{lem3-6} and Markov property, it follows that  for any real number $a\geq
  2 \max\{a(x,y) \; | \; (x,y)\in V(N,\eps)\}$, for any $T\geq 1$ and for any
$(T,\eps,N)$-cluster $\Gamma =
\bigl((t_1,x_1,y_1),\ldots ,(t_n,x_n,y_n)\bigr)$, 
\begin{equation}\label{e10-7}
\frac{1}{a}\log \chi_a(\Gamma) \leq - {\textstyle \sum'
(t_i-t_{i-1})\left( M_N(1, x_i,y_i)  - \varepsilon \right) -
    \sum'' (t_i-t_{i-1}) I^*_0 N/\eps}
\end{equation}
where  $\sum'$ denotes the sum over all $i=1,\ldots,n$ for which
$M_N(1,x_i,y_i) < +\infty$ and $\sum''$ is the sum over all those
$i=1,\ldots,n$ for which $M_N(1,x_i,y_i) = +\infty$.
 
If ${\textstyle \sum'' (t_i-t_{i-1}) \geq T \eps/N,}$ because of  Relation (\ref{e10-7}), $
\log \chi_a(\Gamma)) \leq  - T I^*_0 a + T a \eps$ 
and hence, inequality~(\ref{e9-7}) holds.
 
Furthermore, let us consider for every $1\leq i \leq n$, for which $M_N(1,x_i,y_i) <
+\infty$, a continuous function $\phi_i : [0,1] \to {\cal E}$ with  $\phi_i(0) = x_i$,
$\phi_i(1) = y_1$  and $\sup_t|\phi_(t)| \geq N$, such that   
\begin{equation}\label{e11-7}
I_{[0,1]}(\phi_i) = M_N(1, x_i,y_i).
\end{equation}
Such a function exists because the set of all  continuous functions $\phi : [0,1] \to
{\cal E}$ with $\phi(0) = x_i$, 
$\phi(1) = y_1$  and $\sup_t|\phi(t)| \geq N$ is closed in $D([0,1],{\cal E})$,
and because  the level  sets $\{\phi : \phi(0)=x_i, \; I_{[0,1]}(\phi) \leq c\}$  are
compact (see the definition of a good rate function in Section~\ref{sec4}). 
 
 Let $J(\Gamma)$ be the set of all those $1\leq i\leq n$  for which  the set $\{t : \phi_i(t) 
= 0\}$ is empty and $M_N(1,x_i,y_i) <
+\infty$.  Because of  relation (\ref{e2-7}), $
M_N(1, x_i,y_i) > I^*_0/\eps$ for all  $i\not\in J(\Gamma)$. 
When $\sum'_{i\not\in J(\Gamma)} (t_i-t_{i-1}) \geq T \eps$, using Inequality (\ref{e10-7}) we
get  
$
\log \chi_a(\Gamma) \leq  - T I^*_0 a + T a \eps $
and hence, inequality~(\ref{e9-7}) holds. 
 
To prove our lemma, it is
sufficient now to verify  Inequality~(\ref{e9-7}) in the case when 
\begin{equation}\label{e12-7}
 {\textstyle \sum'_{i\not\in J(\Gamma)} }(t_i-t_{i-1}) < T \eps \quad \; \text{ and } \quad \; 
 {\textstyle \sum'' (t_i-t_{i-1}) < T \eps/N}. 
\end{equation}
For this, we construct an absolutely continuous function $\phi : [0,\tilde{t}]\to{\cal E}$ with
$\phi(0)=\phi(\tilde{t})=0$ for which the Lebesgue measure of the set 
$\{s\in[0,\tilde{t}] : \phi(s) \not= 0\}$ is greater than $T(1 - 2\eps)$ and  
\begin{equation}\label{e13-7}
I_{[0,\tilde{t}]}(\phi) ~\leq~  {\textstyle \sum'} 
(t_i-t_{i-1}) M_N(1, x_i,y_i)  + 2T\eps C/c + 2(bN +T\eps)C/c. 
\end{equation}
 
\medskip
{\it Construction of the function $\phi$} : 
We let $\hat{x}_i = (t_i-t_{i-1}) x_i$ and
$\hat{y}_i = (t_i-t_{i-1}) y_i$ for every $i=1,\ldots,n$  and 
we define 
$0~=~t'_0 < t''_0 < t'_1 < t''_1 < \ldots < t'_n < t''_n =\tilde{t}$ by setting
$t''_0 = \Delta_0$ and $t''_i = t'_i + \Delta_i$ for all for $i=0,\ldots,n$, 
where $
\Delta_0 = ( t_1 \eps + b N)/c$, \, $ 
\Delta_n =
((t_n-t_{n-1})\eps + b N)/c$ and  
$
\Delta_i =
(t_{i+1}-t_{i-1})\eps/c$ for $ i=1,\ldots,n-1$, 
and by letting $t'_i = t''_{i-1} + (t_i - t_{i-1})$ if $M_N(1,x_i,y_i) < +\infty$ and $t'_i =
t''_{i-1} + 2N(t_i - t_{i-1})/c$ otherwise. The function $\phi :[0,\tilde{t}]\to {\cal E}$
is defined then in the following way :  
for $t\in[t''_{i-1},t'_i]$, $i=1,\ldots,n$, we let 
\[
\phi(t) = \begin{cases}
 (t_i-t_{i-1}) \phi_i\bigl((t-t''_{i-1})/(t_i-t_{i-1})\bigr) &\text{ if $M_N(1,x_i,y_i) <
+\infty$}, \\ 
\hat{x}_i + c (\hat{y}_i-\hat{x}_i) (t-t''_{i-1})/(2N (t_i-t_{i-1})) &\text{ if $M_N(1,x_i,y_i) =
+\infty$}.
\end{cases}\]
Then for every $1< i < n$, $\phi(t'_i)  = (t_i-t_{i-1}) y_i = \hat{y}_i$ and $\phi(t''_i)
= (t_{i+1}-t_{i})x_{i+1} = \hat{x}_{i+1}$. For $t\in [t'_i,t''_i]$, $i=0,\ldots,n$, we let 
\[
\phi(t) =  \hat{y}_i  + (\hat{x}_{i+1} - \hat{y}_i)(t-t'_i)/(t''_i-t'_i) 
\]
where $ \hat{y}_0= \hat{x}_{n+1} = 0$. The resulting function $\phi$ is absolutely
continuous on $[0,\tilde{t}]$  and $\phi(0) = \phi(\tilde{t}) = 
0$. Moreover,  by
  construction, $\phi(t)\not=0$ for all $t\in \cup_{i\in J(\Gamma)} [t''_{i-1},t'_i ]$.
  Using Relations (\ref{e12-7}) this implies that 
\begin{align*} 
mes\{t : \phi(t) \not= 0\} &\geq~  {\textstyle \sum'_{i\in J(\Gamma)} }(t_i-t_{i-1}) \\ &\geq~ T
-  {\textstyle \sum'_{i\not\in J(\Gamma)} }(t_i-t_{i-1}) -  
{\textstyle \sum'' (t_i-t_{i-1})} ~\geq~ T - 2T\eps.  
\end{align*} 
The last inequality combined  with Lemma~\ref{lem2-4} and Proposition~\ref{prop1-4} show that 
\begin{equation}\label{e15-7}
I_{[0,\tilde{t}]}(\phi) \geq T (1-2\eps) I^*_0.   
\end{equation}
 
\smallskip
 
{\it Proof of Inequality~(\ref{e13-7})} :  Relations (\ref{e1-4}) and 
(\ref{e11-7}) imply that 
\begin{equation}\label{e16-7}
I_{[t''_{i-1},t'_i]}(\phi) = (t_i-t_{i-1})I_{[0,1]}(\phi_i) = 
(t_i-t_{i-1}) M_N(1,x_i,y_i)
\end{equation}
for all $i=1,\ldots,n$ for which  $M_N(1,x_i,y_i)
< +\infty$.  When $M_N(1,x_i,y_i) = +\infty$,  using Assumption ~(A) and Inequality~(\ref{e8-7}) we obtain 
\begin{equation}\label{e17-7}
I_{[t''_{i-1},t'_i]}(\phi) \leq (t'_i-t''_{i-1}) l_2\bigl((y_i-x_i)c/2N\bigr) \leq 2N C (t_i-t_{i-1})/c. 
\end{equation}
 Because of Assumption~(A),  
\[
I_{[t'_i,t''_i]}(\phi) \leq (t''_i-t'_i) l_2\bigl((\hat{x}_{i+1} -
\hat{y}_i)/\Delta_i\bigr)
\]
 for all $i=0,\ldots,n$.
Moreover, for $1\leq i \leq n-1$, 
\begin{align*}
|\hat{x}_{i+1} - \hat{y}_i| &= |(t_{i+1}-t_i) x_{i+1} - (t_i-t_{i-1}) y_i | \\ &\leq
 (t_i-t_{i-1})\delta(x_i,y_i) +  (t_{i+1}-t_i) \delta(x_{i+1},y_{i+1}) \leq
 (t_{i+1}-t_{i-1}) \eps = c \Delta_i 
\end{align*}
where the first inequality holds because by definition of $(T,\eps,N)$-cluster, the set 
\[(t_i-t_{i-1})B\bigl(y_i,\delta(x_i,y_i)\bigr) \cap (t_{i+1}-t_i)
B\bigl(x_{i+1},\delta(x_{i+1},y_{i+1})\bigr)
\]
is non-empty and the second  inequality follows from Relation~(\ref{e6-7}). Similarly, for
$i=0$,
\[
|\hat{x}_{1} - \hat{y}_0| ~=~ t_1 |x_1| ~\leq~ t_1 \delta(x_1,y_1) + bN ~\leq~ t_1\eps + bN ~=~ c
 \Delta_0 
\]
and for $i=n$, 
\begin{align*}
|\hat{x}_{n+1} - \hat{y}_n| &=~ (t_n-t_{n-1}) |y_n| ~\leq~ (t_n-t_{n-1}) \delta(x_n,y_n) +
 bN \\ &\leq~  (t_n-t_{n-1}) \eps + bN ~=~ c \Delta_n 
\end{align*}
because under the hypotheses of our lemma, the sets  $t_1 B(x_1,
\delta(x_1,y_1) \cap B(0,bN)$ and $(t_n-t_{n-1}) B(y_n,
\delta(x_n,y_n)\cap B(0,bN)$ are non-empty. 
 Using Inequality~(\ref{e8-7}) this implies that  $l_2\bigl((\hat{x}_{i+1}
 - \hat{y}_i)/\Delta_i\bigr) \leq C$ and consequently, 
\[
I_{[t'_i,t''_i]}(\phi) ~\leq~ (t''_i-t'_i) C ~=~ \Delta_i C  
\]
for every $i=0,\ldots,n$. The last inequality, Relations (\ref{e12-7}), 
(\ref{e16-7}), (\ref{e17-7}) 
and  the equality 
\[
\sum_{i=0}^n \Delta_i =  \Bigl( 2 b N +  t_1 \eps +
(t_n-t_{n-1})\eps + \eps \sum_{i=1}^{n-1}
(t_{i+1}-t_{i-1})  \Bigr)/c = 2(b N  +  T\eps)/c
\]
show that  
\begin{align*}
I_{[0,\tilde{t}]}(\phi) &= \sum_{i=1}^n I_{[t''_{i-1},t'_i]}(\phi) +
\sum_{i=0}^n I_{[t'_i,t''_i]}(\phi)  \\ &\leq {\textstyle \sum'}
(t_i-t_{i-1})M_N(1,x_i,y_i)  + 2N C {\textstyle \sum''}
(t_i-t_{i-1})/c + \sum_{i=0}^n\Delta_i C \\
&\leq {\textstyle \sum'} (t_i-t_{i-1})M_N(1,x_i,y_i) +  2TC\eps + 2(b N  +
T\eps) C/c.
\end{align*} 
Relation (\ref{e13-7}) is therefore verified.
 
\medskip
 
We are ready now to complete the proof of Lemma~\ref{lem6-6} :  Relations (\ref{e10-7}),
(\ref{e13-7}) and (\ref{e15-7}) imply that   
\begin{align*}
\frac{1}{a}\log \chi_a(\Gamma) &\leq  - {\textstyle \sum'} 
(t_i-t_{i-1}) \bigl(M_N(1, x_i,y_i) - \eps\bigr) \\ 
 &\leq - I_{[0,\tilde{t}]}(\phi) + 2T\eps C/c + 2(bN +T\eps)C/c + \eps {\textstyle \sum'} 
(t_i-t_{i-1})  \\ 
&\leq  - T(1-2\eps) I^*_0 + 4T\eps C/c + 2bN C/c + T\eps. 
\end{align*}
The last inequality proves Relation~(\ref{e9-7}). 
 \end{proof}
 
\medskip
 
Let us complete now the proof of the inequality $\log r^*_e \leq -I^*_0$.
 
\noindent
Inequalities  (\ref{e1-7}) and (\ref{e7-7})  prove that 
\[
\log r^*_e \leq \limsup_{a\to +\infty} \, \lim_{b\to +\infty} \, 
\limsup_{T\to +\infty} \,  \frac{1}{Ta}\log \sum_\Gamma \chi_{a}(\Gamma).  
\]
The sum is taken here over all
$(T,\eps,N)$-clusters $\Gamma =$  
$((t_1,x_1,y_1),\ldots, (t_n,x_n,y_n))$ for
which the sets $t_1 B(x_1,
\delta(x_1,y_1) \cap B(0,bN)$ and $(t_n-t_{n-1}) B(y_n,
\delta(x_n,y_n)\cap B(0,bN)$ are 
non-empty. Using Lemma~\ref{lem4-6} and Lemma~\ref{lem6-6}  we
obtain therefore  
\begin{align*}
\log r^*_e &\leq~ \limsup_{a\to +\infty} \, \lim_{b\to +\infty} \, 
\limsup_{T\to +\infty} \, \left(
\frac{2}{a} \log(2V(N,\eps)) +
\max_\Gamma \, \frac{1}{T a}\log\chi_{a}(\Gamma)\right)\\ 
&\leq~ \limsup_{a\to +\infty} \, \left(
\frac{2}{a} \log(2V(N,\eps)) - (1-2\eps) I^*_0 + 4\eps C/c + \eps \right)\\ 
&\leq~  - (1-2\eps) I^*_0 + 4\eps C/c + \eps 
\end{align*}
Letting finally $\eps\to 0$ we conclude that $\log r^*_0 \leq -I^*_0$.

\section{Example : applications to Jackson networks}\label{sec7}
 
In this section, we apply our results to Jackson networks.
Let us recall the
definition and some well-known results concerning Jackson network.
This is a network with $d$ queues. For $i=1,\ldots,d$, the arrivals and the service times at
the $i$-th queue are Poisson with parameters $\lambda_i$ and $\mu_i$ respectively. All the
Poisson processes are independent. When the customer finish its service at queue $i$, it goes
to the queue $j$ with probability $p_{ij}$. The residual quantity
$
p_{i0} = 1 -  p_{i1}+\cdots +p_{id}
$
is the probability that the customer leaves definitely the network.
We assume that $p_{ii}=0$ for all $i=1,\ldots,d$.
 
Let $X_i(t)$  denote the length of the queue i at time $t$, $i=1,\ldots,d$. Then
$X(t)=(X_1(t),\ldots, X_d(t))$ is a continuous time Markov chain on
$\Z_+^d$ with generator
\[
{\cal G}f(x) =
\sum_{y\in\Z_+^d} q(x,y) (f(y)-f(x)), \quad x\in\Z_+^d,
\]
where
\begin{equation}\label{e1-J}
q(x,y) = \begin{cases}\lambda_i &\mbox{ if } y-x=\epsilon_i, \; \;
  i\in\{1,\ldots,d\},\\
\mu_ip_{i0} &\mbox{ if } y-x= - \epsilon_i, \; \; i\in\{1,\ldots,d\},\\
\mu_ip_{ij} &\mbox{ if } y-x= \epsilon_j-\epsilon_i, \; \; i,j\in\{1,\ldots,d\}, i\not=j,\\
0 &\mbox{ otherwise,}
\end{cases}
\end{equation}
$\epsilon_i$ is the $i$-th unit vector. We set $p_{00}=1$ and $p_{0i}=0$ for all $i\not=
0$, the matrix $(p_{ij}; i,j=0,\ldots,d)$ is then stochastic.
 
\noindent 
{\bf Assumption (J)\; }{\em We  suppose that
\begin{itemize}
\item[1)] the spectral radius of the matrix $(p_{ij}; \, i,j=1,\ldots,d)$ is strictly
less than unity, 
\item[2)] for any $1\leq i\leq d$, there exists $n\in\N$ and $1\leq j\leq d$ such that
$\lambda_j p_{ji}^{(n)} > 0$  where
$p_{ji}^{(n)}$ denotes  the $n$-time transition probability of a Markov chain with
$d+1$ states associated to the stochastic matrix $(p_{ij}; \, i,j=0,\ldots,d)$. 
\end{itemize}}
\medskip
 
The Markov process $(X(t))$ is then irreducible. The system of traffic equations
\[
\nu_i = \lambda_i + \nu_1p_{1i} + \cdots +  \nu_dp_{di} 
\]
has the unique solution $(\nu_i)$~(see~\cite{Kelly}).
The Markov process $(X(t))$ is recurrent if and only if $\nu_i \leq \mu_i$  for all
$i=1,\ldots,d$,  and it is ergodic (positive recurrent) if and only if $\nu_i < \mu_i$ for all 
$i=1,\ldots,d$.
 
Recall that {\it spectral radius} of the process
 $(X(t))$ is defined by
\[
r^* = \inf\left\{ r>0 : \int_0^\infty r^{-t} \P_x(X(t)=y)  \, dt < +\infty, \;
\forall x,y\in\Z^d_+ \right\}
\]
When the process $(X(t))$ is recurrent we have obviously $r^*=1$. For
a transient Markov process, spectral radius $r^*$ shows how fast the process goes
to infinity.
 
{\it Essential spectral radius} $r^*_e$ is defined as the infimum over all those $r>0$ for which there
is a finite set $K\subset\Z^d_+$ such that 
\[
\int_0^\infty r^{-t} \P_x\bigl( X(t)=y, \;\tau (K) > t \bigr)  \, dt < +\infty \quad 
\text{ for all } \quad  x,y\in\Z^d_+ \setminus K
\]
where $\tau(K)$ denotes the first time when the process $(X(t))$ hits the set
$K$. 
 
When the Markov process  $(X(t))$ is recurrent the quantity $r^*_e$ is the infimum
over all those $r>0$ for which there is a finite set $K\subset\Z^d_+$ such that 
\[
\int_0^\infty r^{-t} \P_x\bigl(\tau (K) > t \bigr)  \, dt < +\infty \quad 
\text{ for all } \quad  x\in\Z^d_+ \setminus K.
\]
The hypotheses of Proposition~\ref{pr3-3}  are satisfied  here because the Markov 
process $(X(t))$ has uniformly bounded jumps. 
 
The Markov process $(X(t))$ satisfies the hypotheses~(A) : 
the family of scaled Markov processes $(Z_a(t), \; t\in[0,T])\stackrel{\text{def.}}{=}
(X(at)/a, \; t\in[0,T])$ satisfies sample path large deviation  
principle (see~\cite{A-D,D-E, Ignatiouk:04}) with a good rate function  having  an
integral form :  
\[
I_{[0,T]}(\phi) =
\int_0^T L(\phi(t),\dot{\phi}(t)) \, dt
\]
for every absolutely continuous function $\phi :[0,T]\to\R_+^d$, and $I_{[0,T]}(\phi)
=+\infty$ otherwise. 
The local rate function $L(x,v)$ can be represented in several ways,
see~\cite{A-D,Ignatiouk:01}:   for
$\Lambda\subset\{1,\ldots,d\}$ and $x=(x_1,\ldots, x_d)\in\R^d_+$ with $x_i>0$
for $i\in\Lambda$ and $x_i=0$ for $i\not\in\Lambda$,
\[
L(x,v) ~=~ L_\Lambda(v) ~=~ \sup_{\alpha\in{\cal B}_\Lambda} \bigl(\langle \alpha, v\rangle -
R(\alpha) \bigr) ~=~ \sup_{\alpha\in\R^d} \left(\langle \alpha, v\rangle -
\max_{\Lambda': \Lambda\subset\Lambda'} R_{\Lambda'}(\alpha) \right)
\]
for all $v\in\R^d$. $\langle \cdot, \cdot \rangle$ denotes here the usual scalar product in $\R^d$,
\[
R(\alpha) ~=~ \sum_{j=1}^d \mu_i\Bigl( \sum_{j\not=i} p_{ij}e^{\alpha_j-\alpha_i} +
p_{i0}e^{-\alpha_i} - 1\Bigr) + \sum_{j=1}^d \lambda_i(e^{\alpha_i}-1),
\]
\[
R_{\Lambda'}(\alpha) ~=~ \sum_{j\in\Lambda'} \mu_i\Bigl( \sum_{j\not=i} p_{ij}e^{\alpha_j-\alpha_i} +
p_{i0}e^{-\alpha_i} - 1\Bigr) + \sum_{j=1}^d \lambda_i(e^{\alpha_i}-1)
\]
and ${\cal B}_\Lambda$ is the set of those all  $\alpha = (\alpha_1, \ldots,\alpha_d) \in\R^d$
for which  
\[
\alpha_i\leq \log\Bigl( \, \sum_{j=1}^d p_{ij}e^{\alpha_j} +
p_{i0} \Bigr) \quad \text{ for all $ i\not\in\Lambda$.}
\]
The hypotheses~$(a_3)$ are satisfied with  
\[
l_1(v) = \sup_{\alpha\in{\cal B}_\emptyset} \bigl(\langle \alpha, v\rangle -
R(\alpha) \bigr) = \sup_{\alpha\in\R^d}  \left(\langle \alpha, v\rangle -
\max_{\Lambda\subset\{1,\ldots,d\}} R_\Lambda(\alpha) \right) 
\]
and 
\[
l_2(v) = R^*(\alpha) = \sup_{\alpha\in\R^d}  \left(\langle \alpha, v\rangle -
R(\alpha) \right).
\] 
Under Assumption~(J), the convex conjugate $R^*$ of the function $R$ is finite everywhere
on $\R^d$ (see \cite{Ignatiouk:01}, Lemma~10.1), and for every $\delta > 0$, 
\[
l_1(v) ~\geq~ \sup_{|\alpha|\leq \delta}  \langle \alpha, v\rangle - \sup_{|\alpha|\leq
  \delta} \; \max_{\Lambda\subset\{1,\ldots,d\}} R_\Lambda(\alpha) ~\geq~ \delta |v| - \sup_{|\alpha|\leq
  \delta} \; \max_{\Lambda\subset\{1,\ldots,d\}} R_\Lambda(\alpha).
\]
which implies that 
\[
\lim_n \inf_{|v|\geq n} l_1(v)/|v|  > 0. 
\]
 
Using Theorem~\ref{th1} and
Theorem~\ref{th2} it follows that 
\begin{cor}\label{cor1-7}
Under the hypotheses (J),  
\begin{equation}\label{e2-J}
\log r^* ~=~ -  L_\emptyset(0) ~=~ - \inf_{\phi}
I_{[0,1]}(\phi)     
\end{equation}
where the infimum is taken over all continuous functions $\phi : [0,1]\to \R_+^d$ with
$\phi(0)=\phi(1)$ and 
\begin{equation}\label{e3-J}
\log r^*_e ~=~ - \inf_{\phi~:~
  \phi(t)\not=0, \, \forall 0<t<1} I_{[0,1]}(\phi) 
\end{equation}
where the infimum is taken over all continuous functions $\phi : [0,1]\to \R_+^d$ with
$\phi(0)=\phi(1)$ such that $\phi(t)\not=0$ for all $0<t<1$.
\end{cor}
 
Using this result we calculate explicitly the quantity $r^*_e$ for $d=1$ and for
$d=2$. This is a subject of the following propositions. 
 
\begin{prop}\label{prop1-J} For $d=1$, $\log r^*_e = \inf_{\alpha\in\R} R(\alpha)$. 
\end{prop}
\begin{proof} Indeed, for any continuous function $\phi : [0,1]\to \R_+$ with
$\phi(0)=\phi(1)$ and $\phi(t)\not=0$ for all $0<t<1$, 
\[
I_{[0,1]}(\phi) = \int_0^1 R^*(\dot\phi(t))  \, dt ~\geq~ R^*(0) = -\inf_{\alpha\in\R} R(\alpha) 
\]
because the convex conjugate $R^*$ of the function $R$  is convex. For a constant function
$\phi(t)\equiv x$ with $x>0$,  
\[
I_{[0,1]}(\phi) = R^*(0) = -\inf_{\alpha\in\R} R(\alpha).
\]
This proves that the right hand side of  Relation~(\ref{e3-J})  equals 
$\inf_{\alpha\in\R} R(\alpha)$ and hence, $\log r^*_e = \inf_{\alpha\in\R} R(\alpha)$. 
\end{proof}

\begin{prop}\label{prop2-J} For $d=2$,
\begin{equation}\label{e4-J}
\log r^*_e ~=~ - \inf_i L_{\{i\}}(0) ~=~   - (1- p_{12}p_{21}) \min\{ (\sqrt{\mu_1} - \sqrt{\nu_1})^2,
(\sqrt{\mu_2} - \sqrt{\nu_2})^2\} 
\end{equation}
 if the Markov process $(X(t))$ is ergodic, and 
\begin{equation}\label{e5-J}
\log r^*_e  ~=~ - \inf_i L_{\{i\}}(0) ~=~ \log r^*  ~=~ - L_\emptyset(0) 
\end{equation}
otherwise.
\end{prop}
 
To prove Proposition~\ref{prop2-J} we consider the following lemmas.
 
\begin{lemma}\label{lem1-J}
For all
$\alpha =(\alpha_1,\alpha_2)\in{\cal B}_{\{1\}}$ and
$\beta=(\beta_1,\beta_2)\in{\cal B}_{\{2\}}$ such that $\alpha_1\geq\beta_1$ and
$\beta_2\geq \alpha_2$, the following inequality holds
\begin{equation}\label{e7-J}
\log r^*_e \leq \max\bigl\{R(\alpha),R(\beta)\bigr\}.
\end{equation}
\end{lemma}
\begin{proof} The function $R(\cdot)$ being continuous,
it is sufficient to prove our lemma for the case when $\alpha_1>\beta_1$ and
$\alpha_2<\beta_2$.
 
Recall that ${\cal B}_{\{1\}}$  is the set of all $(\alpha_1,\alpha_2)\in\R^2$
with $\alpha_2 \leq \log(p_{21}e^{\alpha_1}+p_{20})$, and
${\cal B}_{\{2\}}$ is the set of all $(\alpha_1,\alpha_2)\in\R^2$
with $\alpha_1 \leq \log(p_{12}e^{\alpha_2}+p_{10})$.
For given $\alpha =(\alpha_1,\alpha_2)\in{\cal B}_{\{1\}}$ and
$\beta=(\beta_1,\beta_2)\in{\cal B}_{\{2\}}$ satisfying the inequalities $\alpha_1 > \beta_1$  and
$\beta_2 > \alpha_2$, let $
f(x) = \exp{\langle \alpha,x \rangle } + \exp{\langle \beta , x \rangle}$.
Then for $x=(x_1,x_2)\in\Z_+^2$ with $x_1>0$ and $ x_2>0,$ we have
\begin{equation}\label{e8-J}
{\cal G}f(x) \, = \, R(\alpha) \exp{\langle \alpha,x \rangle } +
R(\beta)\exp{\langle \beta,x \rangle } \, \leq \, \max\{R(\alpha),R(\beta)\}  f(x).
\end{equation}
Furthermore, for $\theta = (\theta_1,\theta_2)\in\R^2$,  denote
$c_1(\theta) = -\mu_2\left(p_{21}e^{\theta_1-\theta_2} + p_{20}e^{-\theta_2} -
1 \right)$. The inequality $\alpha_2\leq \log(p_{21}e^{\alpha_1}+p_{20})$
implies that $c_1(\alpha) \leq 0$ and therefore, for $x_1>0$ and $ x_2=0,$ we obtain
\begin{align*}
{\cal G}f(x) &= (R(\alpha) + c_1(\alpha)) \exp(\alpha_1x_1) +
(R(\beta) + c_1(\beta))\exp(\beta_1x_1)\\
&\leq R(\alpha) \exp(\alpha_1x_1) +
(R(\beta) + c_1(\beta))\exp(\beta_1x_1)
\end{align*}
and
\[
{\cal G}f(x)/f(x) \leq R(\alpha) + (R(\beta) + c_1(\beta))\exp(\beta_1x_1 - \alpha_1x_1).
\]
The right hand side of the last inequality tends to $R(\alpha)$ as
$x_1\to\infty$ because $\alpha_1>\beta_1$, and hence, for any $\eps > 0$ there is $N_1(\eps)>0$ such that
for $x_1 >  N_1(\eps)$ and $x_2=0$,
\begin{equation}\label{e9-J}
{\cal G}f(x) \leq \bigl( R(\alpha) + \eps \bigr) f(x).
\end{equation}
The same arguments show that for any $\eps > 0$ there is $N_2(\eps)>0$ such that for
$x_2 >  N_1(\eps)$ and
$x_1=0$,
\[
{\cal G}f(x) \leq \bigl( R(\beta) + \eps \bigr) f(x).
\]
The last inequality combined with  Relations (\ref{e8-J}) and  (\ref{e9-J})  shows that 
\[
{\cal G}f(x) \leq \Bigl( \max\{R(\alpha),R(\beta)\} +
\varepsilon\Bigr) f(x)
\]
for all $x=(x_1,x_2)\in\Z_+^2$ with 
$x_1+x_2 > N(\eps) = \max\{N_1(\eps),N_1(\eps)\}$. This implies that for all $t> 0$ and
$x=(x_1,x_2)\in\Z^2_+$ with  
$x_1+x_2 > \max\{N_1(\eps),N_1(\eps)\}$, 
\[
\E_x( f(X(t)) ; \; \tau(W) > 0 ) \leq \exp \Bigl( t \max\{R(\alpha),R(\beta)\} +
t\eps\Bigr)  f(x)
\]
where $\tau(W)$ denotes the first time when the process $(X(t))$ hits the set
$W = \{(x_1,x_2)\in\Z^2_+: x_1+x_2\leq N(\eps)\}$. Using Proposition~\ref{pr1-3}
we conclude that 
\[
\log r^*_e \leq \max\{R(\alpha),R(\beta)\} + \eps
\]
and letting $\eps\to 0$  we obtain inequality (\ref{e7-J}).
\end{proof}
 
\medskip
 
Let ${\cal R}$ denote
 the infimum of $\max\{ R(\alpha_1,\alpha_2), R(\beta_1,\beta_2)\}$  over all
$(\alpha_1,\alpha_2)\in{\cal B}_{\{1\}}$ and $(\beta_1,\beta_2)\in{\cal B}_{\{2\}}$
with $\alpha_1\geq\beta_1$ and
$\beta_2\geq \alpha_2$. 
\begin{lemma}\label{lem2-J} When the Markov process $(X(t))$
is ergodic,
\begin{equation}\label{e10-J}
{\cal R} \leq -\min_i L_{\{i\}}(0) =  - (1- p_{12}p_{21}) \min\{ (\sqrt{\mu_1} - \sqrt{\nu_1})^2,
(\sqrt{\mu_2} - \sqrt{\nu_2})^2\}.
\end{equation}
\end{lemma}
\begin{proof}
Recall that the Markov process $(X(t))$ is ergodic if and only if $\nu_1 < \mu_1$ and
$\nu_2 < \mu_2$ where $(\nu_1,\nu_2)$ is a unique solution of the traffic equations $
\nu_1 = \lambda_1 + \nu_2p_{21}$ and $\nu_2 = \lambda_2 + \nu_1p_{12}$.
 
Let $\partial {\cal B}_{\{i\}}$ denote the boundary of the set ${\cal B}_{\{i\}}$.
Straightforward
calculations show that for $(\alpha_1,\alpha_2)\in\partial {\cal B}_{\{i\}}$,
\[
R(\alpha_1,\alpha_2) =  (1- p_{12}p_{21}) \bigl(\mu_i(e^{-\alpha_i}-1) +
\nu_i(e^{\alpha_i}-1)\bigr), \quad i=1,2.
\]
The minimum of the function $R(\cdot)$ on the boundary $\partial {\cal B}_{\{1\}}$ is achieved
at the point $\alpha^*=(\alpha^*_1,\alpha^*_2) $ with
\[
\alpha^*_1 = \log\sqrt{\mu_1/\nu_1} > 0 \quad  \alpha^*_2=
\log\bigr(p_{21}\sqrt{\mu_1/\nu_1}+p_{20} \bigl) \geq 0,
\]
and equals
\[
R(\alpha^*) = -  (1- p_{12}p_{21}) (\sqrt{\mu_1} - \sqrt{\nu_1})^2.
\]
The minimum of the function $R(\cdot)$ on the boundary $\partial {\cal B}_{\{2\}}$ is achieved
at the point $\beta^*=(\beta^*_1,\beta^*_2) $ with
\[
\beta^*_2 = \log\sqrt{\mu_2/\nu_2}  > 0, \quad \beta^*_1=
\log\bigr(p_{12}\sqrt{\mu_2/\nu_2}+p_{10} \bigl) \geq 0,
\]
and equals
\[
R(\beta^*) = -  (1- p_{12}p_{21}) (\sqrt{\mu_2} - \sqrt{\nu_2})^2.
\]
Without any restriction of generality we will suppose that $ R(\alpha^*)\geq R(\beta^*)$
which implies that
\[
 L_2(0)
\geq -R(\beta^*) \geq - R(\alpha^*) =  (1- p_{12}p_{21}) \min\{(\sqrt{\mu_1} - \sqrt{\nu_1})^2,
(\sqrt{\mu_2} - \sqrt{\nu_2})^2\}.
\]
Hence, to prove the equality
\begin{equation}\label{e11-J}
\min_i L_{\{i\}}(0) =  (1- p_{12}p_{21}) \min\{(\sqrt{\mu_1} - \sqrt{\nu_1})^2,
(\sqrt{\mu_2} - \sqrt{\nu_2})^2\}
\end{equation}
it is sufficient to show that $L_{\{1\}}(0) = - R(\alpha^*)$. For this, we consider the set
$ S = \{ \alpha=(\alpha_1,\alpha_2) : R(\alpha) < R(\alpha^*)\}$ and we notice that
$\beta^*\not\in{\cal B}_{\{1\}}$ because $\beta^*\in{\cal B}_{\{2\}}$ and  
$\beta^*_1>0$, $\beta^*_2>0$
while for all $\beta = (\beta_1,\beta_2)\in{\cal B}_{\{1\}}\cap{\cal B}_{\{2\}}$,
$\beta_1\leq 0$ and $\beta_2 \leq 0$. 
Since $\alpha^*\in\partial{\cal B}_{\{1\}}$ this
implies that $t\alpha^* + (1-t)\beta^* \not\in {\cal B}_{\{1\}}$
for all $0<t<1$ because  the set $\R^2\setminus{\cal B}_{\{1\}}$ is convex. The function
$R(\cdot)$ being strictly convex, the inequality 
$R(\alpha^*)\geq R(\beta^*)$ implies that $t\alpha^* + (1-t)\beta^* \in S$
for all $0< t < 1$. The set $S$ being connected we conclude that $S\cap {\cal B}_{\{1\}} =
\emptyset$ because the boundary of the set ${\cal B}_{\{1\}}$ has no intersection with
$S$. This proves that the point $\alpha^*$ achieves the infimum of the function $R(\cdot)$ on
${\cal B}_{\{1\}}$ and consequently, $L_{\{1\}}(0) = - R(\alpha^*)$. Relation (\ref{e11-J}) is therefore 
proved.
 
\smallskip
 
Notice furthermore that 
\[
\beta^*_2 \, = \, \sqrt{\mu_2/\nu_2} \; \geq \; 1 + (\sqrt{\mu_1/\nu_1} - 1) \sqrt{\nu_1/\nu_2}
\;
\geq \; p_{21} \sqrt{\mu_1/\nu_1}  + p_{20} \, = \, \alpha^*_2
\]
where the first inequality follows from the inequality $ R(\alpha^*)\geq R(\beta^*)$ and the second
inequality holds because $\nu_1 = \lambda_1 + \nu_2 p_{21} \geq \nu_2 p_{21}^2$.
When $\alpha^*_1 \geq \beta^*_1$ we obtain therefore
${\cal R} \leq \max\{R(\alpha^*), R(\beta^*)\} = R(\alpha^*) = - L_{\{1\}}(0)$ and hence,
relation  (\ref{e10-J}) holds.
\begin{figure}[ht]
\resizebox{7cm}{6cm}{\includegraphics{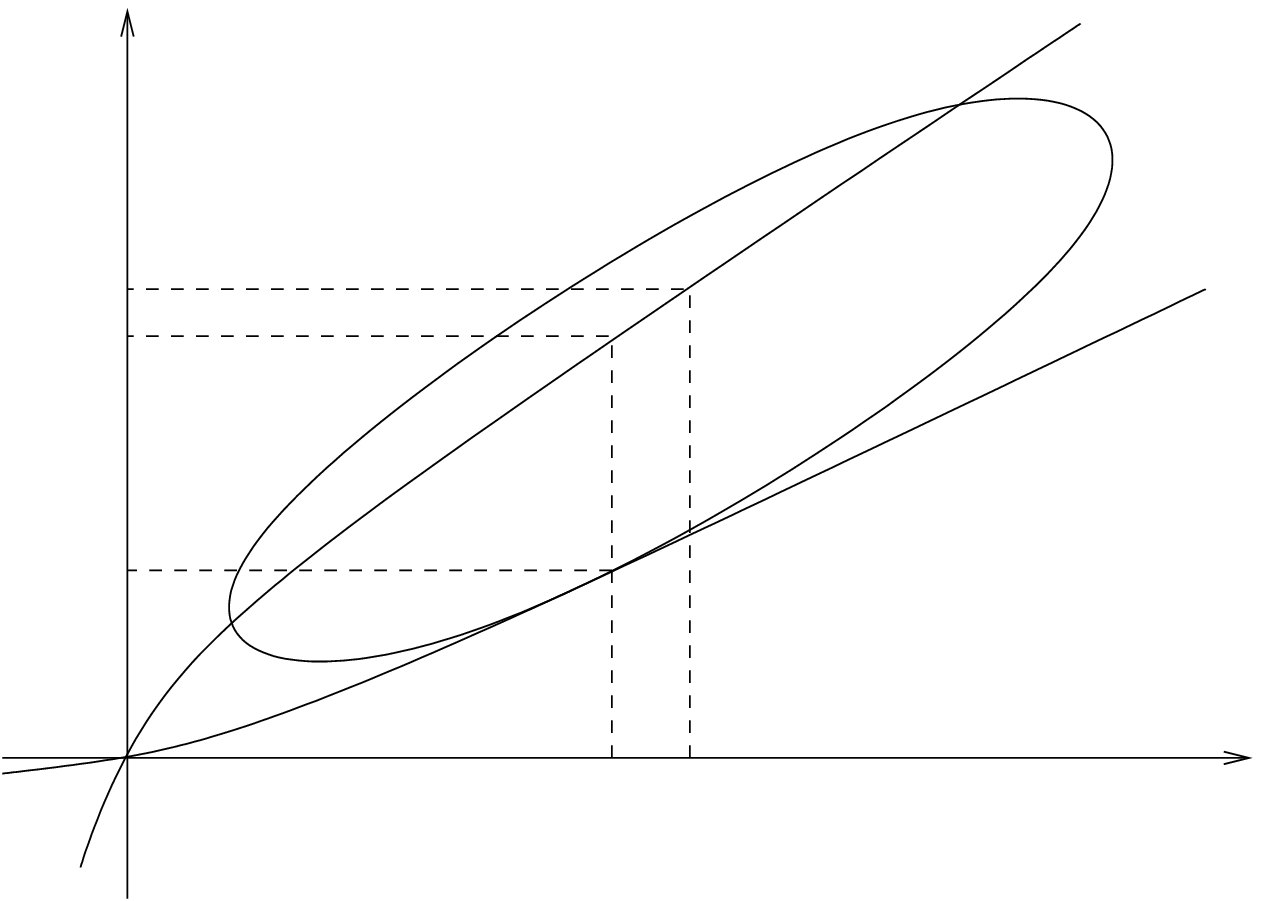}}
\put(-110,20){$\alpha^*_1$}
\put(-195,100){$\hat\beta_2$}
\put(-195,118){$\beta^*_2$}
\put(-195,62){$\alpha^*_2$}
\put(-90,20){$\beta^*_1$}
\put(-82,100){$S$}
\put(-50,60){${\cal B}_1$}
\put(-220,20){$\partial{\cal B}_1$}
\put(-200,-5){$\partial{\cal B}_2$}
\put(-140,125){${\cal B}_2$}
\caption{ }
\end{figure}
Suppose now that $\alpha^*_1 < \beta^*_1$ (see Figure~1), then $p_{12}\not=0$ because
otherwise, $\beta^*_1 = 0 < \alpha^*_1$. For
$(\hat\beta_1,\hat\beta_2)\in\partial{\cal B}_{\{2\}}$ with $\hat\beta_1 = \alpha^*_1$ and
$\hat\beta_2 = \log(e^{\alpha^*_1} - p_{10}) -\log p_{12}$,
we have $\hat\beta_2\geq  \alpha^*_1 \geq \alpha^*_2$ which implies that  ${\cal R}\leq \max\{R(\alpha^*),
R(\hat\beta)\}$. Moreover, using the inequality $
\hat\beta_2\leq \log(e^{\beta^*_1} - p_{10}) -\log p_{12} = \beta^*_2$ we obtain
\begin{align*}
R(\hat\beta) &=   (1- p_{12}p_{21}) 
(\nu_2 - \mu_2  e^{-\hat\beta_2} )(e^{\hat\beta_2}-1) \\ &\leq  (1- p_{12}p_{21})  (\nu_2  - \mu_2
e^{-\beta^*_2}) (e^{\alpha^*_1} - 1)/p_{12} \\ 
&\leq -  (1- p_{12}p_{21}) (\sqrt{\mu_2/\nu_2}-1) (\sqrt{\mu_1/\nu_1} - 1)\, \nu_2 /p_{12} \\ &\leq -
 (1- p_{12}p_{21}) (\sqrt{\mu_1/\nu_1} - 1)^2 \, = \, R(\alpha^*) 
\end{align*}
where the last inequality holds because $\sqrt{\mu_2/\nu_2} - 1 \geq
(\sqrt{\mu_1/\nu_1} - 1)\sqrt{\nu_1/\nu_2} > 0$ and
$\sqrt{\nu_2} \geq \sqrt{\nu_1p_{12}} \geq \sqrt{\nu_1}p_{12}.$ This proves that
${\cal R}\leq R(\alpha^*) = - L_{\{1\}}(0)$ and consequently,  Relation (\ref{e10-J}) holds.
\end{proof}
 
\noindent 
{\bf Proof of Proposition~\ref{prop2-J}}.   
Relation~(\ref{e2-J})  imply that
\begin{equation}\label{e12-J}
\log r^* = -L_\emptyset(0) = - \inf_{\phi(0)=\phi(1)} I_{[0,1]}(\phi)
\end{equation}
where the infimum is taken over all continuous functions $\phi :[0,1]\to \R^2_+$
with $\phi(0)=\phi(1)$. Because of Relation~(\ref{e3-J}), 
\begin{equation}\label{e13-J}
\log r^*_e = - \inf_{\phi} I_{[0,1]}(\phi) 
\end{equation}
where the infimum is taken over all those continuous functions $\phi :[0,1]\to \R^2_+$
for which $\phi(0)=\phi(1)$ and $\phi(t)\not= 0$ for all $0<t<1$.
 
When the Markov process $(X(t))$ is non-ergodic,  Theorem~2 of
Ignatiouk~\cite{Ignatiouk:01} proves that 
\[
L_\emptyset(0) = \min_i L_{\{i\}}(0) = \min\{  I_{[0,1]}(\phi), I_{[0,1]}(\psi)\} 
\]
where $\phi(t)=x'=(1,0)\in\R_+^2$ and $\psi(t)=x'' = (0,1)\in\R_+^2$ for all $t\in[0,1]$.  
Using Relations (\ref{e12-J}) and (\ref{e13-J}) this implies that 
\[
\log r^* = - L_\emptyset(0) = -\min_i L_{\{i\}}(0) = - \min\{  I_{[0,1]}(\phi), I_{[0,1]}(\psi)\}  \, \leq
  \,
\log r^*_e \leq \log r^*.
\]
and consequently, Relation (\ref{e5-J}) holds.  
 
Suppose now that the Markov process $(X(t))$ is ergodic. Then using Lemma~\ref{lem1-J}
and Lemma~\ref{lem2-J} it follows that 
\begin{align*}
\log r^*_e \leq {\cal R} \leq -\min_i L_{\{i\}}(0) &=  -  (1- p_{12}p_{21}) \min\{ (\sqrt{\mu_1} - 
\sqrt{\nu_1})^2,
(\sqrt{\mu_2} - \sqrt{\nu_2})^2\} \\ & = - \min\{  I_{[0,1]}(\phi), I_{[0,1]}(\psi)\}.   
\end{align*}
The last relation combined with Inequality (\ref{e13-J}) proves Relation (\ref{e4-J}).

\providecommand{\bysame}{\leavevmode\hbox to3em{\hrulefill}\thinspace}
\providecommand{\MR}{\relax\ifhmode\unskip\space\fi MR }
\providecommand{\MRhref}[2]{%
  \href{http://www.ams.org/mathscinet-getitem?mr=#1}{#2}
}
\providecommand{\href}[2]{#2}


\begin{thebibliography}{10}
 
\bibitem{A-D}
Rami Atar and Paul Dupuis, \emph{Large deviations and queueing networks:
  methods for rate function identification}, Stochastic Processes and their
  Applications \textbf{84} (1999), no.~2, 255--296. \MR{1 719 274}
 
\bibitem{Billingsley}
Patrick Billingsley, \emph{Convergence of probability measures}, Wiley series
  in probability and mathematical statistics, John Wiley \& Sons Ltd, New York,
  1968.
 
\bibitem{Delcoigne:01}
Franck Delcoigne and Arnaud~de La~Fortelle, \emph{Large deviations for polling
  systems}, Tech. Report RR-3892, {INRIA}, March 2000, {\em
  http://www.inria.fr/RRRT/RR-3892.html}.
 
\bibitem{D-Z}
Amir Dembo and Ofer Zeitouni, \emph{Large deviations techniques and
  applications}, Springer-Verlag, New York, 1998.
 
\bibitem{D-E}
Paul Dupuis and Richard~S. Ellis, \emph{The large deviation principle for a
  general class of queueing systems. {I}}, Transactions of the American
  Mathematical Society \textbf{347} (1995), no.~8, 2689--2751.
 
\bibitem{Grillo:01}
G.~Grillo, \emph{On {P}ersson's theorem in local {D}irichlet spaces},
  Zeitschrift f\"ur Analysis und ihre Anwendungen. Journal for Analysis and its
  Applications \textbf{17} (1998), no.~2, 329--338.
 
\bibitem{Ignatiouk:01}
Irina Ignatiouk-Robert, \emph{The large deviations of {J}ackson networks},
  Annals of Applied Probability \textbf{10} (2000), no.~3, 962--1001.
 
\bibitem{Ignatiouk:04}
\bysame, \emph{Large deviations for processes with discontinuous statistics},
  Tech. Report 14/04, {Universit\'e de Cergy-Pontoise, UMR 8088}, March 2004.
 
\bibitem{Kelly}
F.~P. Kelly, \emph{Reversibility and stochastic networks}, Wiley, Chichester,
  1979.
 
\bibitem{Malyshev-Minlos}
V.~A. Malyshev and R.~A. Minlos, \emph{Gibbs random fields}, Kluwer Academic
  Publishers Group, Dordrecht, 1991, Cluster expansions, Translated from the
  Russian by R. Koteck\'y and P. Holick\'y.
 
\bibitem{Malyshev-Spieksma}
V.~A. Malyshev and F.M. Spieksma, \emph{Intrinsic convergence rate of countable
  {M}arkov chains}, Markov Processes and Related Fields \textbf{1} (1995),
  203--266.
 
\bibitem{Nummelin1}
Esa Nummelin, \emph{General irreducible {M}arkov chains and nonnegative
  operators}, Cambridge University Press, Cambridge, 1984.
 
\bibitem{Rivasseau:01}
Vincent Rivasseau, \emph{From perturbative to constructive renormalization},
  Princeton University Press, Princeton, NJ, 1991.
 
\bibitem{Seneta}
E.~Seneta, \emph{Nonnegative matrices and {M}arkov chains}, second ed.,
  Springer-Verlag, New York, 1981. \MR{85i:60058}
 
\bibitem{S-W}
A.~Shwartz and A.~Weiss, \emph{Large deviations for performance analysis},
  Stochastic Modeling Series, Chapman \& Hall, London, London, 1995.
 
\bibitem{Stroock}
D.W. Stroock, \emph{On the spectrum of {M}arkov semigroups and the existence of
  invariant measures}, Functional Analysis in {M}arkov processes (M.~Fukushima,
  ed.), Springer Verlag, 1981, pp.~287--307.
 
\bibitem{Vere-Jones:1}
D.~Vere-Jones, \emph{Ergodic properties of nonnegative matrices-{I}}, Pacific
  Journal of mathematic \textbf{22} (1967), no.~2, 361--386.
 
\bibitem{Vere-Jones:2}
\bysame, \emph{Ergodic properties of nonnegative matrices-{II}}, Pacific
  Journal of mathematic \textbf{26} (1968), no.~3, 601--620.
 
\bibitem{Varadhan-Williams}
R.J. Williams and S.R.S. Varadhan, \emph{Brownian motion in a wedge with
  oblique reflection}, Comm. Pure Appl. Math. (1985), no.~24, 147--225.
 
\bibitem{Woess}
Wolfgang Woess, \emph{Random walks on infinite graphs and groups}, Cambridge
  University Press, Cambridge, 2000.
 
\bibitem{LimingWu}
Liming Wu, \emph{Essential spectral radius for markov semigroups (i): discrete
  time case}, Probab. Theory Relat. Fields \textbf{128} (2004), 255--321.
 
\end{thebibliography}
\end{document}